\documentclass[a4paper, 11pt]{article}

\usepackage[utf8]{inputenc} % Encodage des caractères
\usepackage[T1]{fontenc}    % Encodage des fontes
\usepackage[english]{babel} % Langue du document
\usepackage{amsmath, amssymb, amsthm} % Packages AMS pour les mathématiques
\usepackage{geometry}       % Pour régler les marges
\usepackage{graphicx}       % Pour inclure des images
\usepackage{hyperref}       % Pour les liens hypertextes et les références
\usepackage{tocloft}        % Pour personnaliser le sommaire
\usepackage{braket}
\usepackage{cleveref}
\usepackage{tikz}
\usepackage{nicematrix}
\usepackage[numbers]{natbib}

\newcommand{\bk}[2]{\braket{#1|#2}}

\newcommand{\ov}[4]{o_{#1 #2 #3 #4}^v}
\newcommand{\ou}[4]{o_{#1 #2 #3 #4}^u}
\newcommand{\ow}[4]{o_{#1 #2 #3 #4}^w}

\numberwithin{equation}{section}

\title{Eigenvector Overlaps of Random Covariance Matrices and their Submatrices}
\author{Attal Elie\thanks{Ecole Polytechnique, CMAP, elie.attal@polytechnique.edu}, Allez Romain\thanks{Qube Research \& Technologies, romain.allez@qube-rt.com}}
\date{\today}

\begin{document}

% Page de titre
\maketitle

\begin{abstract}
    We consider the singular vectors of any $m \times n$ submatrix of a rectangular $M \times N$ Gaussian matrix and study their asymptotic overlaps with those of the full matrix, in the macroscopic regime where $N \,/\, M\,$, $m \,/\, M$ as well as $n \,/\, N$ converge to fixed ratios. Our method makes use of the dynamics of the singular vectors and of specific resolvents when the matrix coefficients follow Brownian trajectories. We obtain explicit forms for the limiting rescaled mean squared overlaps for right and left singular vectors in the bulk of both spectra, for any initial matrix $A\,$. When it is null, this corresponds to the Marchenko-Pastur setup for covariance matrices, and our formulas simplify into Cauchy-like functions.
\end{abstract}

\section{Introduction}
\label{sec:Introduction}

Suppose $A$ is a deterministic $M \times N$ matrix with $M \geq N$ and $B_t$ has the same dimensions and contains independent Brownian motions. The matrix
\begin{equation}
    \label{eq:X_def}
X_t := A + \frac{1}{\sqrt{N}} \, B_t
\end{equation}
can be viewed as a noisy observation of $A\,$. For $m < M$ and $n < N\,$, we are interested in comparing $X_t$ with $\tilde{X}_t\,$, defined by
\begin{equation}
\label{eq:X_tilde_definition}
\tilde{X}_t^{ij} = 
\begin{cases}
X_t^{ij}\,,\,\,\, &\textit{if $i \leq m$ and $j \leq n\,$,} \\
0\,\hspace{4mm},\,\,\, &\textit{otherwise,}
\end{cases}
\end{equation}
through their singular vectors. We focus on the case $n \leq m$ so that $\tilde{X}_t$ has rank $n$ almost surely.

\begin{equation*}
    \tilde{X}_t = 
    \begin{pmatrix}
        X_t^{11} & \cdots & X_t^{1n} & 0 & \cdots & 0 \\
        \vdots &  & \vdots & \vdots &  & \vdots \\
        \vdots & \ddots & \vdots & \vdots & \ddots & \vdots \\
        \vdots &  & \vdots & \vdots &  & \vdots \\
        X_t^{m1} & \cdots & X_t^{mn} & 0 & \cdots & 0 \\
        0 & \cdots & 0 & 0 & \cdots & 0 \\
        \vdots & \ddots  & \vdots & \vdots & \ddots & \vdots \\
        0 & \cdots & 0 & 0 & \cdots & 0
    \end{pmatrix} \,.
\end{equation*}
\newpage
\noindent
Let us introduce:
\begin{itemize}
\item[$\bullet$] $X_t = U_t\,\Lambda_t \,V_t^T$ the Singular Values Decomposition (SVD) of $X_t\,$, with singular values $\sqrt{\lambda_1^t} \geq ... \geq \sqrt{\lambda_N^t} > 0\,$, left singular vectors $u_1^t\,, ...\,, u_M^t$ and right singular vectors $v_1^t\,, ...\,, v_N^t\,$.
\item[$\bullet$] $\tilde{X}_t = \tilde{U}_t \, M_t \, \tilde{V}_t^T$ the SVD of $\tilde{X}_t\,$, with singular values $\sqrt{\mu_1^t} \geq ... \geq \sqrt{\mu_n^t} > 0 = \sqrt{\mu_{n+1}^t} = ... = \sqrt{\mu_N^t}\,$, left singular vectors $\tilde{u}_1^t\,, ... \,,\tilde{u}_M^t$ and right singular vectors $\tilde{v}_1^t\,, ... \,,\tilde{v}_N^t\,$. We add the following condition: the null space of $\tilde{X}_t^T$ can be divided into two parts. The singular vectors $\tilde{u}_{n+1}^t\,, ... \,,\tilde{u}_m^t$ have all their $M-m$ last components equal to zero, representing the fact that $m \geq n$ so that the first $n$ columns do not form a free family of vectors. Furthermore, the vectors $\tilde{u}_{m+1}^t\,, ...\,, \tilde{u}_M^t$ have all their first $m$ components equal to zero, representing the part of the null space due to the shape of $\tilde{X}_t$ and its $M - m$ null columns. Specifically, the latter can be seen as $e_{m+1}\,, ...\,, e_{M}\,$ (where $e_i$ has all his coefficients null except the i-th which equals 1). Note that this condition corresponds to taking certain linear combinations of the vectors of the null space, and therefore does not modify the formulas obtained for the singular vectors associated with non-zero singular values.
\end{itemize}
We are interested in the limiting behaviour of the overlaps $\bk{\tilde{u}_i^t}{u_j^t}$ and $\bk{\tilde{v}_i^t}{v_j^t}$ for any $t\,$, as $M\,, N \,, m\,, n \to \infty$ with $N\,/\,M \to q\,$, as well as $n \,/\, N \to \alpha$ and $m \,/\,M \to \beta\,,$ i.e. the macroscopic regime. Specifically, we study these limits for singular vectors in the bulk of both spectra. This is equivalent to studying the overlaps between the eigenvectors of the square matrices $R_t := X_t^T\, X_t\,$, $\tilde{R}_t := \tilde{X}_t^T\, \tilde{X}_t\,$, $L_t := X_t\, X_t^T$ and $\tilde{L}_t := \tilde{X}_t \, \tilde{X}_t^T\,$ which are empirical covariance matrices of $X_t$ or $\tilde{X}_t\,$. When $A$ is null, one can view $X_t$ as a dataset of $M$ independent samples of $N$ independent Gaussian variables of mean zero and variance $t\,$, and $\tilde{X}_t$ is a subselection of a macroscopic number of samples and features. Our work allows one to compare the Principal Component Analysis (PCA) of $X_t$ with $\tilde{X}_t$'s eigenvector by eigenvector, under the assumption of independent features, which corresponds to the Marchenko-Pastur setup. Note that similarly to \cite{attal2024interlacing}, the time $t$ is the variance of the noise added to $A\,$, but it is also a way to derive dynamics that allow us to obtain our results.

As mentioned in \cite{attal2024interlacing}, there is no trivial deterministic relation between the eigenvectors of a symmetric matrix and those of one of its principal minors. In that context, the Random Matrix Theory approach has proved to be a powerful tool allowing to obtain explicit asymptotic formulas for the expectations of the squared overlaps. The case of Wishart matrices we are considering here is no different, we expect to obtain similar results for the overlaps of left and right singular vectors using random matrices.

Moreover, the use of random matrices accounts for the noise measured on top of a relevant signal. The Marchenko-Pastur distribution of singular values (see \citep{marchenko1967distribution,tao2012topics,potters2020first}) for perturbed data or image such as (\ref{eq:X_def}) has been widely used for denoising in many different contexts, including MRI images \citep{veraart2016denoising, veraart2016diffusion, zhu2022denoise}, financial data \citep{bouchaud2009financial, laloux2000random, utsugi2004random} and wireless communications \citep{bai2010spectral, tulino2004random}. Other results that focus on the eigenvalues of Wishart matrices (squared singular values of $X_t$ when $A \equiv 0$ in our setup) such as the BBP phase transition \cite{baik2005phase} and the Tracy-Widom law for extreme eigenvalues \cite{johnstone2001distribution} have found applications in various fields \citep{patterson2006population, majumdar2010extreme}. Although these results mainly focus on the eigenvalues of such random matrices, their eigenvectors have gained interest over the years. The main focus is to derive estimators of the population covariance matrix while observing a sample covariance, such as in \citep{ledoit2011eigenvectors, mestre2008improved, bun2016rotational, lin2024eigenvector}. Additionally, minors of Wishart matrices have been increasingly studied in recent years, with applications in conditional independence in covariance matrices \cite{drton2008moments}, compressed sensing \citep{cai2021asymptotic, jiang2024largest} and percolation theory \cite{adler2013random, dieker2008largest}.

In our previous work \cite{attal2024interlacing}, we derived explicit formulas in the context of symmetric Gaussian matrices for the limiting rescaled mean squared overlaps between the eigenvectors of a principal submatrix and those of the full matrix. Our approach was based on analysing the eigenvector flow under the Dyson Brownian motion and deriving the dynamics of a specific resolvent. However, these findings were confined to symmetric matrices and their principal minors. In the present article, we extend this method to the singular vectors of rectangular Gaussian matrices, or equivalently, to the eigenvectors of Wishart matrices. By examining the singular vectors' dynamics in this context, we establish analogous results for the limiting overlaps in the macroscopic regime. 

Our work therefore reaches two main domains of application. On the one hand, we study the information contained in a subimage of a rectangular noisy image through their singular vectors. On the other hand, we establish a link between the Principal Component Analysis (PCA) of a sample covariance matrix with identity population covariance, and the PCA obtained when removing a macroscopic number of features or samples. In particular, we believe our results can bring new insights into Incremental PCA algorithms \citep{hall1998incremental, artac2002incremental, weng2003candid}, PCA with missing data \citep{nelson1996missing, folch2015pca}, Risk Management or Portfolio Optimization by financial sector \citep{dhingra2024sectoral} or time-dependent PCA methods \citep{li2000recursive, de2015overview}.

In Section \ref{sec:Dynamics}, we introduce the dynamics of the eigenvalues and the eigenvectors and derive the correlation structure of the different Brownian motions in presence. We then recall some results on the Stieltjes transforms of the spectral densities and their limiting Burgers equations. The special case $A \equiv 0$ is shown to be that of the Marchenko-Pastur distribution. Section \ref{sec:Overlaps} contains the resolution of our problem. We introduce the quantities we want to study, and define three resolvents that have forms similar to the one used in \cite{attal2024interlacing}. We prove that in the scaling limit, they become solutions of a deterministic system of coupled differential equations (\ref{eq:PDE_system}) that we are able to solve explicitly. Using an inversion formula, we obtain explicit forms for the limiting rescaled mean squared overlaps, for a general matrix $A$ (see (\ref{eq:General_solution})). In the case $A \equiv 0\,$, we have the following limits (\ref{eq:result_A_null}) for $\lambda\,,\mu > 0\,$,
\setlength{\jot}{20pt}
\begin{align*}
&N \, \mathbb{E}\left[\bk{\tilde{v}_{i_N}^t}{v_{j_N}^t}^2\right] \,\longrightarrow\, q \, \frac{(1 - \alpha) \, t \,\bar{\mu} + \alpha \, (1 - \beta)\, t \, \bar{\lambda} + (1 - \alpha \beta)\, (\alpha + \frac{1}{q})\, t^2}{(1 - \alpha \beta)^2 \, t^2 + q \, (\bar{\lambda} - \bar{\mu})\, (\alpha \beta \, \bar{\lambda} - \bar{\mu})}\,,\\
&N \, \mathbb{E}\left[\bk{\tilde{u}_{i_N}^t}{u_{j_N}^t}^2 \right]\,\longrightarrow\, q\,\frac{(1-\beta)\, t \,\bar{\mu} + \beta \, (1 - \alpha)\, t \, \bar{\lambda} + (1 - \alpha \beta)\, (1 + \frac{\beta}{q})\, t^2}{(1 - \alpha \beta)^2\, t^2 + q\,(\bar{\lambda} - \bar{\mu})\, (\alpha \beta \, \bar{\lambda} - \bar{\mu})}\,,\\
&N \, \mathbb{E}\left[\bk{\tilde{v}_{i_N}^t}{v_{j_N}^t}\, \bk{\tilde{u}_{i_N}^t}{u_{j_N}^t} \right] \,\longrightarrow\, q\,\frac{(1 - \alpha \beta)\, t \, \sqrt{\lambda\, \mu}}{(1 - \alpha \beta)^2\, t^2 + q\,(\bar{\lambda} - \bar{\mu})\, (\alpha \beta \, \bar{\lambda} - \bar{\mu})}\,,
\end{align*}
as $M\,,N\,,m\,,n \to \infty$ with $\left(\frac{N}{M}\,,\frac{n}{N}\,,\frac{m}{M}\right) \to \left(q \,, \alpha \,, \beta\right)$ as well as $\mu_{i_N}^t \to \mu\,$, $\lambda_{j_N}^t \to \lambda$ and using the notations $\bar{\mu} := \mu - \left(\alpha + \frac{\beta}{q}\right)\, t$ and $\bar{\lambda} := \lambda - \left(1 + \frac{1}{q}\right)\, t\,$.

\section{Eigenvalue and Eigenvector Dynamics}
\label{sec:Dynamics}

In 1989, Bru (\citep{bru1989diffusions}) derived the dynamics of the eigenvalues and eigenvectors of $R_t = X_t^T \, X_t\,$. For any $1 \leq j \leq N\,$, we have
\begin{equation}
\label{eq:Bru_dynamics_eigenvalues}
d\lambda_j^t = \frac{2}{\sqrt{N}} \sqrt{\lambda_j^t}\, db_j(t) + \frac{M}{N}\, dt + \frac{1}{N}\, \sum_{\substack{k = 1 \\ k \neq j}}^N \,\frac{\lambda_j^t + \lambda_k^t}{\lambda_j^t - \lambda_k^t}\,dt \,,
\end{equation}
\begin{equation}
\label{eq:Bru_dynamics_eigenvectors}
dv_j^t = -\frac{1}{2N} \, \sum_{\substack{k = 1 \\ k \neq j}}^N \, \frac{\lambda_j^t + \lambda_k^t}{(\lambda_j^t - \lambda_k^t)^2}\,v_j^t\,dt + \frac{1}{\sqrt{N}} \, \sum_{\substack{k =1 \\ k \neq j}}^N \, \frac{\sqrt{\lambda_j^t} \, dw_{jk}(t) + \sqrt{\lambda_k^t}\, dw_{kj}(t)}{\lambda_j^t - \lambda_k^t} \, v_k^t \,,
\end{equation}
where $\left\{b_j\, |\, 1 \leq j \leq N\right\}$ and $\left\{w_{jk} \, |\, 1 \leq j \leq M\,, 1 \leq k \leq N\,, j \neq k\right\}$ are two families of independent Brownian motions, independent of each other. Specifically, in the proof of these dynamics given by Bru, we can identify these processes as $db_j(t) = \bk{u_j^t}{dB_t \,v_j^t}$ and $dw_{jk}(t) = \bk{u_j^t}{dB_t \,v_k^t}\,$. These dynamics are different from those obtained in the symmetric Brownian case, i.e. the Dyson Brownian motion \citep{dyson1962brownian, potters2020first, tao2012topics, grabiner1999brownian}, but we can find some similarities. First, we remark that the eigenvalues are still subject to a repulsion force, which is not exactly inversely proportional to their distance. Moreover, since the family of Brownian motions $dw$ is independent of $db\,$, the eigenvectors' dynamics can be seen as diffusion processes in a random environment, given by the eigenvalues trajectories, which is also the case for the Dyson Brownian motion. Note that this is due to the fact that the Brownian motions in $B_t$ are uncorrelated, otherwise, the coefficients of the population covariance matrix in the $v^t$ and $u^t$ bases would appear in both dynamics.

By replacing $X_t$ with $X_t^T\,$, we can obtain the dynamics of the left singular vectors $u_j^t$ for $1\leq j \leq M\,$. They are distinguished into two cases depending on whether $ j \leq N$ or whether $j > N$ (corresponding to a vector in the null space of $X_t^T$). For $1 \leq j \leq N\,$, we have
\begin{align*}
du_j^t &= -\frac{1}{2N} \, \sum_{\substack{k =1 \\ k \neq j}}^N \, \frac{\lambda_j^t + \lambda_k^t}{(\lambda_j^t - \lambda_k^t)^2} \, u_j^t \, dt + \frac{N - M}{2N \lambda_j^t} \, u_j^t \, dt \\
&\quad + \frac{1}{\sqrt{N}} \, \sum_{\substack{k = 1 \\ k \neq j}}^N \, \frac{\sqrt{\lambda_j^t}\, dw_{kj}(t) + \sqrt{\lambda_k^t}\,dw_{jk}(t)}{\lambda_j^t - \lambda_k^t} \, u_k^t 
 + \frac{1}{\sqrt{N}} \, \sum_{k = N+1}^M \, \frac{dw_{kj}(t)}{\sqrt{\lambda_j^t}}\, u_k^t \,,
\end{align*}
whereas for $N+1 \leq j \leq M\,$,
\begin{equation*}
du_j^t = -\frac{1}{2N} \, \sum_{\substack{k = 1}}^N \frac{1}{\lambda_k^t} \, u_j^t \, dt + \frac{1}{\sqrt{N}} \, \sum_{k=1}^N \,\frac{dw_{jk}(t)}{\sqrt{\lambda_k^t}}\,u_k^t \,.
\end{equation*}
Note that the roles of $dw_{jk}$ and $dw_{kj}$ are exchanged for the left singular vectors. Additionally, these dynamics are identical to
\begin{equation*}
du_j^t = - \frac{1}{2N} \, \sum_{\substack{k = 1 \\ k \neq j}}^M \, \frac{\lambda_j^t + \lambda_k^t}{(\lambda_j^t - \lambda_k^t)^2} \, u_j^t \, dt + \frac{1}{\sqrt{N}}\, \sum_{\substack{k = 1 \\ k \neq j}}^M \, \frac{\sqrt{\lambda_j^t}\, dw_{kj}(t) + \sqrt{\lambda_k^t} \, dw_{jk}(t)}{\lambda_j^t - \lambda_k^t} \, u_k^t \,,
\end{equation*}
for any $1 \leq j \leq M\,$, if we set $\lambda_{N+1}^t = ... = \lambda_M^t = 0$ and using the convention $0\,/\, 0 = 0\,$. This form will be used throughout this paper to simplify our computations. Obviously, the notation $w_{jk}$ with $k > N$ is not properly defined, but with our convention it is always multiplied by a null factor.

The truncated matrix $\tilde{R}_t = \tilde{X}_t^T\, \tilde{X}_t$ has null coefficients outside of its $n \times n$ top left submatrix, and has rank $n$ almost surely. Therefore, its eigenvectors associated with non-zero eigenvalues only have non-zero coefficients on their first $n$ components, and inversely, its eigenvectors in the null space have all their first $n$ components equal to zero. Consequently, there is no interaction with the null space and we can deal with $\tilde{v}_1^t\,, ...\,, \tilde{v}_n^t$ only. They behave the same way the $v_i^t$ do if we replace $N$ with $n$ and $M$ with $m$ (the scaling remains in $1 \,/\,\sqrt{N}$). Similarly, the dynamics of the $\mu_i^t$ can be derived from those of the  $\lambda_j^t\,$, meaning we have for any $1\leq i \leq n\,$,
\begin{equation*}
d\mu_i^t = \frac{2}{\sqrt{N}} \sqrt{\mu_i^t} \, d\tilde{b}_i(t) + \frac{m}{N} \, dt + \frac{1}{N} \, \sum_{\substack{l = 1 \\ l \neq i}}^n \, \frac{\mu_i^t + \mu_l^t}{\mu_i^t - \mu_l^t}\,dt \,,
\end{equation*}
\begin{equation*}
d\tilde{v}_i^t = -\frac{1}{2N} \, \sum_{\substack{l = 1 \\ l \neq i}}^n \, \frac{\mu_i^t + \mu_l^t}{(\mu_i^t - \mu_l^t)^2} \, \tilde{v}_i^t \, dt + \frac{1}{\sqrt{N}} \, \sum_{\substack{l = 1 \\ l \neq i }}^n \, \frac{\sqrt{\mu_i^t} \, d\tilde{w}_{il}(t) + \sqrt{\mu_l^t} \, d\tilde{w}_{li}(t)}{\mu_i^t - \mu_l^t} \, \tilde{v}_l^t \,,
\end{equation*}
where $\left\{d\tilde{b}_i := \bk{\tilde{u}_i^t}{d\tilde{B}_t \, \tilde{v}_i^t} \,|\,\, 1 \leq i \leq n \right\}$ and $\left\{d\tilde{w}_{il} := \bk{\tilde{u}_i^t}{d\tilde{B}_t \, \tilde{v}_l^t} \,|\,\, 1 \leq i \leq m \,, 1 \leq l \leq n \,, i \neq l\right\}$ are independent of each other. Here we defined by $\tilde{B}_t$ the truncated version of the Brownian matrix $B_t\,$, the way we defined $\tilde{X}_t$ from $X_t\,$. For $\tilde{u}_i^t\,$, using the convention $\mu_{n+1}^t = ... = \mu_{m}^t = 0$ and $0\,/\,0 = 0\,$, we get for any $1 \leq i \leq m\,$,
\begin{equation*}
    d\tilde{u}_i^t = - \frac{1}{2N}\, \sum_{\substack{l = 1 \\ l \neq i}}^{m}\, \frac{\mu_i^t + \mu_l^t}{(\mu_i^t - \mu_l^t)^2}\, \tilde{u}_i^t \, dt + \frac{1}{\sqrt{N}}\, \sum_{\substack{l=1 \\ l \neq i}}^{m}\, \frac{\sqrt{\mu_i^t}\, d\tilde{w}_{li}(t) + \sqrt{\mu_l^t}\, d\tilde{w}_{il}(t)}{\mu_i^t - \mu_l^t}\, \tilde{u}_l^t \,.
\end{equation*}

Finally, if we want to study the overlaps between the singular vectors $X_t$ and those of $\tilde{X}_t\,$, we need to compute the correlations between the different Brownian motions. In Appendix \ref{subsec:App_correlation}, we prove that
\begin{equation*}
\bk{u_j^t}{dB_t \, v_k^t} \, \bk{\tilde{u}_i^t}{d\tilde{B}_t \, \tilde{v}_l^t} = \bk{\tilde{u}_i^t}{u_j^t} \, \bk{\tilde{v}_l^t}{v_k^t} \, dt \,,
\end{equation*}
for any $(i,l,j,k) \in \{1\,; ...\,; m\} \times \{1\,;...\,;n\} \times \{1\,;...\,;M\} \times \{1\,;...\,;N\}\,$. Thus, when the following correlations are properly defined, we have
\begin{align*}
\begin{cases}
dw_{jk}(t) \, d\tilde{b}_i(t) = \bk{\tilde{u}_i^t}{u_j^t} \,  \bk{\tilde{v}_i^t}{v_k^t} \, dt \,,\\
db_j(t) \, d\tilde{w}_{il}(t) = \bk{\tilde{u}_i^t}{u_j^t} \, \bk{\tilde{v}_l^t}{v_j^t} \, dt \,,\\
db_j(t) \, d\tilde{b}_i(t) = \bk{\tilde{u}_i^t}{u_j^t} \, \bk{\tilde{v}_i^t}{v_j^t} \, dt \,,\\
dw_{jk}(t) \, d\tilde{w}_{il}(t) = \bk{\tilde{u}_i^t}{u_j^t} \, \bk{\tilde{v}_l^t}{v_k^t} \, dt\,.
\end{cases}
\end{align*}

Since our work focuses on eigenvectors in the bulk, we need to make the following assumption: the spectrum of $A^T A\,$, i.e. $\lambda_1^0 \geq ... \geq \lambda_N^0\,$ (recall that $X_0 = A$), has an empirical distribution converging to a continuous density $\rho(\cdot, 0)\,$:
\begin{equation*}
\frac{1}{N} \sum_{j=1}^N \, \delta_{\lambda_j^0}(d\lambda) \longrightarrow \rho(\lambda,0) \, d\lambda \,.
\end{equation*}
Similarly, we assume
\begin{equation*}
\frac{1}{n} \, \sum_{i=1}^n \, \delta_{\mu_i^0}(d\lambda) \longrightarrow \tilde{\rho}(\mu,0) \, d\mu \,,
\end{equation*}
where $\tilde{\rho}(\cdot,0)$ is also continuous. For any time $t\,$, we denote the (continuous) limiting density of $R_t$'s spectrum (respectively of the non-zero part of $\tilde{R}_t$'s spectrum) by $\rho(\cdot,t)$ (respectively $\tilde{\rho}(\cdot,t)$). We can therefore define the Stieltjes transforms associated with both spectra,
$$
G_N(z,t) := \frac{1}{N}\, \sum_{j=1}^N \, \frac{1}{z - \lambda_j^t} \quad \text{and} \quad \tilde{G}_N(\tilde{z},t) := \frac{1}{n} \sum_{i = 1}^n \, \frac{1}{\tilde{z} - \mu_i^t}\,,
$$
and write their respective limits as $N \to \infty$ as
$$
G(z,t) := \int_{\mathbb{R}} \, \frac{\rho(\lambda,t)}{z - \lambda} \, d\lambda \quad \text{and} \quad \tilde{G}(\tilde{z},t) := \int_{\mathbb{R}} \, \frac{\tilde{\rho}(\mu,t)}{\tilde{z} - \mu} \, d\mu \,.
$$
These functions, defined for $z\,, \tilde{z} \in \mathbb{C} \setminus \mathbb{R}\,$, are classical tools used to study the limiting behaviour of the spectral densities. Indeed, one can recover $\rho$ from $G$ using the Sokhotski-Plemelj formula
\begin{equation}
    \label{eq:Stieltjes_inversion_formula}
\lim_{\varepsilon \to 0^+} \, G(\lambda \pm i \, \varepsilon\,, t) = v(\lambda, t) \mp i \pi \, \rho(\lambda, t)\,,
\end{equation}
where $v(\lambda, t) :=  P.V. \,\int_{\mathbb{R}} \, \frac{\rho(\lambda', t)}{\lambda - \lambda'}\, d\lambda'$ is the Hilbert transform of $\rho\,$ and $P.V.$ denotes Cauchy's principal value.

Applying Itô's lemma, we find, in the scaling limit, the following Burgers equation (see Appendix \ref{subsubsec:App_deriving_burgers}),
\begin{equation}
\label{eq:G_Burgers}
\partial_tG = \left(1 - \frac{1}{q} - 2z \, G \right) \, \partial_z G - G^2 \,,
\end{equation}
which was originally found in \cite{duvillard2001large}\,. Notice that this limiting differential equation is deterministic, which confirms the intuition that the spectral density becomes deterministic in the scaling limit and the eigenvalues stick to their quantiles. In the context of symmetric Gaussian matrices of \cite{attal2024interlacing}, we also find a Burgers equation, however it does not contain the additive $G^2$ term. Using the method of characteristics (see Appendix \ref{subsubsec:App_solving_burgers}), equation (\ref{eq:G_Burgers}) can be solved, leading to an implicit equation on $G$ in the general case:
\begin{equation}
\label{eq:G_implicit}
G(z,t) = \frac{G\left(z_t\, z_t'\,, 0 \right)}{1 + t \, G\left(z_t\, z_t'\,, 0 \right)} \,,
\end{equation}
where $z_t := 1 - t\, G(z,t)$ and $z_t' := z \, \left(1 - t \, G(z,t)\right) - \left(q^{-1} - 1\right)\, t\,$. In the case $A \equiv 0\,$, we have $G(z, 0) = 1 \,/\, z$ and the equation gives $G(z,t)$ as a zero of a second-order polynomial. We can find the correct root due to the fact that $G(z,t) \sim 1\,/\,z$ as $|z| \to \infty\,$. We obtain
\begin{equation*}
G(z,t) = \frac{z - (\frac{1}{q} - 1)\,t - \sqrt{(z - (1 + \frac{1}{\sqrt{q}})^2\,t)\,(z - (1 - \frac{1}{\sqrt{q}})^2\,t)}}{2 z t} \,.
\end{equation*}
It corresponds to the Stieltjes transform of the Marchenko-Pastur distribution
\begin{equation*}
\rho(\lambda,t) = \frac{\sqrt{(\lambda_+ - \lambda)(\lambda- \lambda_-)}}{2\pi \lambda t}
\end{equation*}
with $ \lambda_{\pm} = (1 \pm 1 \, / \, \sqrt{q})^2\, t\,$, see \cite{marchenko1967distribution}\,.

Similarly, we find for $\tilde{G}$ the limiting equation
\begin{equation}
\label{eq:G_tilde_Burgers}
\partial_t \tilde{G} = \left(\alpha - \frac{\beta}{q} - 2 \alpha \tilde{z} \, \tilde{G} \right) \, \partial_{\tilde{z}}\tilde{G} - \alpha\, \tilde{G}^2 \,,
\end{equation}
leading to the implicit equation
\begin{equation}
\label{eq:G_tilde_implicit}
\tilde{G}(z,t) = \frac{\tilde{G}\left(\tilde{z}_t\, \tilde{z}_t'\,,0   \right)}{1 + \alpha t\, \tilde{G}\left(\tilde{z}_t\, \tilde{z}_t'\,,0   \right)} \,,
\end{equation}
where $\tilde{z}_t := 1 - \alpha t \, \tilde{G}\left(\tilde{z}, t\right)$ and $\tilde{z}_t' := \tilde{z}\, \left(1 - \alpha t \, \tilde{G}\left(\tilde{z}, t\right)\right) - \left(\beta \,/\,q - \alpha\right)\, t\,$.
When $A \equiv 0\,$, it is the Stieltjes transform of the Marchenko-Pastur density with $\mu_{\pm} = \left(1 \pm \sqrt{\beta \,/\, \alpha q}\right)^2 \alpha t\,$.

These equations will be pivotal for the remainder of our computations.

\section{Limiting Behaviour of the Overlaps}
\label{sec:Overlaps}
\subsection{The General Case}
We now define the quantities under investigation. Let us introduce the notations
\begin{align*}
&V_{ij}(t) := \bk{\tilde{v}_i^t}{v_j^t}^2 \,, \quad \text{for $1 \leq i \leq n$ and $1 \leq j \leq N\,$,} \\
&U_{ij}(t) := \bk{\tilde{u}_i^t}{u_j^t}^2 \,, \quad \text{for $1 \leq i \leq m$ and $1 \leq j \leq M\,$,} \\
&W_{ij}(t) := \bk{\tilde{v}_i^t}{v_j^t} \, \bk{\tilde{u}_i^t}{u_j^t} \,, \quad \text{for $1 \leq i \leq n$ and $1 \leq j \leq N\,$.}
\end{align*}
The normalisation constraints of the orthonormal bases indicate that these objects vanish as $1 \,/\, N$ in the bulk, so that our goal is to compute the limits of $N \, \mathbb{E}\left[V_{ij}(t)\right]\,$, $N \, \mathbb{E}\left[U_{ij}(t)\right]\,$ and $N \, \mathbb{E}\left[W_{ij}(t)\right]\,$. More precisely, if $i_n \, / \, n \to x \in [0\,,1]$ and $j_N \,/\,N \to y \in [0\,,1]\,$, we have $N \, \mathbb{E}\left[V_{i_n\,j_N}(t)\right] \to V(x,y,t)\,$, where the limiting overlapping function $V$ is the object we want to explicit. Similarly, $N \, \mathbb{E}\left[U_{i_n \, j_N}(t)\right] \to U(x,y,t)$ and $N \, \mathbb{E}\left[W_{i_n \, j_N}(t)\right] \to W(x,y,t)\,$. For the left singular vectors, we have three other cases:
\begin{itemize}
\item[$\bullet$] If $n+1 \leq i_n \leq m$ and $j_N\,/\, N \to y \in [0\,,1]\,$, then $N \, \mathbb{E}\left[U_{i_n\,j_N}(t)\right] \to U^{(1)}(y,t)$ which does not depend on $i$ because the roles of $u_{n+1}^t\,, ...\,, u_m^t$ can be exchanged.
\item[$\bullet$] If $i_n\,/\,n \to x \in [0\,,1]$ and $N+1 \leq j_N \leq M\,$, then $N\, \mathbb{E}\left[U_{i_n\,j_N}(t)\right] \to U^{(2)}(x,t)\,$.
\item[$\bullet$] If $n+1 \leq i_n \leq m$ and $N+1 \leq j_N \leq M\,$, then $N\,\mathbb{E}\left[U_{i_n\,j_N}(t)\right] \to U^{(3)}(t)\,$.
\end{itemize}
We can define the quantile functions $\lambda(\cdot, t)$ and $\mu(\cdot, t)$ of the limiting spectral densities at time $t$ as 
\begin{equation*}
    \label{eq:quantile_def}
    x = \int_{\lambda(x,t)}^{\infty} \, \rho(\lambda,t)\, d\lambda = \int_{\mu(x,t)}^{\infty} \, \tilde{\rho}(\mu,t)\, d\mu\,.
\end{equation*}
They allow us to define more suitable target functions using a change of variable. We define $\bar{V}\,$, $\bar{U}$ and $\bar{W}$ with $\bar{V}\left(\mu(x,t)\,, \lambda(y,t)\,, t\right) = V(x\,, y\,, t)\,$. The function $\bar{U}$ can be extended to the three other cases with:
\begin{itemize}
    \item[$\bullet$] $\bar{U}\left(0\,, \lambda(y,t)\,, t\right) = U^{(1)}(y,t)\,,$
    \item[$\bullet$] $\bar{U}\left(\mu(x,t)\,, 0 \,, t\right) = U^{(2)}(x,t)\,,$
    \item[$\bullet$] $\bar{U}(0\,,0\,,t) = U^{(3)}(t)\,.$  
\end{itemize}
Note that in most applications, we are only interested in the overlaps of singular vectors associated with non-zero singular values. Moreover, numerical simulations of overlaps involving singular vectors of the null space can vary depending on the chosen vector, but their roles are theoretically exchangeable in the scaling limit. We include these cases in our study (only for the left singular vectors, since we treat the case $M \geq N$ and $m \geq n$) because contrary to the symmetric case of \cite{attal2024interlacing}, the vectors of the null space appear in the dynamics of Section \ref{sec:Dynamics}, and are therefore needed to achieve the calculations.

Similarly to \cite{attal2024interlacing}\,, the dynamics of $V_{ij}\,$, $U_{ij}$ and $W_{ij}$ (see Appendix \ref{subsec:App_ito_overlaps}) are difficult to deal with directly. In order to find an explicit expression for their limits, we introduce three complex functions of the variables $z\,,\tilde{z} \in \mathbb{C} \setminus \mathbb{R}\,$,
$$
S_V^{(N)}(z,\tilde{z},t) := \frac{1}{N} \, \sum_{i=1}^n\, \sum_{j=1}^N \, \frac{V_{ij}(t)}{(\tilde{z} - \mu_i^t)(z-\lambda_j^t)} \,,
$$
$$
S_U^{(N)}(z,\tilde{z},t) := \frac{1}{N} \, \sum_{i=1}^m\, \sum_{j=1}^M \, \frac{U_{ij}(t)}{(\tilde{z} - \mu_i^t)(z-\lambda_j^t)} \,,
$$
$$
S_W^{(N)}(z,\tilde{z},t) := \frac{1}{N} \, \sum_{i=1}^n\, \sum_{j=1}^N \, \frac{\sqrt{\mu_i^t \lambda_j^t} \, W_{ij}(t)}{(\tilde{z} - \mu_i^t)(z-\lambda_j^t)} \,.
$$
They have self-averaging properties as the sum of many different random variables, and still encode all the information of the squared overlaps as Stieltjes transforms. They play a role similar to the resolvents used in \cite{attal2024interlacing}, \cite{bun2018overlaps}, \cite{ledoit2011eigenvectors} and \cite{pacco2023overlaps}\,. We typically expect these quantities to converge to deterministic integrals involving the goal functions $\bar{V}\,, \bar{U}$ and $\bar{W}\,$.

This intuition is confirmed in Appendix \ref{subsec:deriving_system}, as we show that applying Itô's lemma to all three resolvent gives a deterministic system of coupled partial differential equations in the scaling limit:
\begin{equation}
    \label{eq:PDE_system}
    \begin{cases}
    \partial_t S_V = g(z,t)\, \partial_z S_V + \tilde{g}(\tilde{z}, t) \, \partial_{\tilde{z}}S_V + \left(2\, S_W - G(z,t) - \alpha \, \tilde{G}(\tilde{z},t)\right)\, S_V \\
    \\
    \partial_t S_U = g(z,t)\, \partial_z S_U + \tilde{g}(\tilde{z}, t) \, \partial_{\tilde{z}}S_U + \left(2\, S_W - \frac{\frac{1}{q} - 1}{z} - G(z,t) - \frac{\frac{\beta}{q} - \alpha}{\tilde{z}} -  \alpha \, \tilde{G}(\tilde{z},t)\right)\, S_U \\
    \\
    \partial_t S_W = g(z,t)\, \partial_z S_W + \tilde{g}(\tilde{z}, t) \, \partial_{\tilde{z}}S_W + S_W^2 + z \tilde{z}\,\,S_V\, S_U\,,
    \end{cases}
\end{equation}
where $g(z,t) := 1 - \frac{1}{q} - 2z\, G(z,t)$ and $\tilde{g}(\tilde{z},t) := \alpha - \frac{\beta}{q} - 2\alpha \tilde{z}\,\tilde{G}(\tilde{z},t)\,$. Since the characteristics of these equations are the same as those of (\ref{eq:G_Burgers}) and (\ref{eq:G_tilde_Burgers}), we can solve them using the method of characteristics. If we introduce the notation $S^0(t) := S(z_t\, z_t'\,, \tilde{z}_t \, \tilde{z}_t'\,, 0)$ for any of our three Stieltjes transforms, then our solutions are given by 
$$
\begin{cases}
    S_V(z, \tilde{z}, t) = \frac{z_t\, \tilde{z}_t\, S_V^0(t)}{\left(1 - t \, S_W^0(t)\right)^2 - z_t\,z_t'\,\tilde{z}_t\, \tilde{z}_t'\, S_V^0(t)\, S_U^0(t) \, t^2} \\
    \\
    S_U(z, \tilde{z}, t) = \frac{z_t'\, \tilde{z}_t'\, S_U^0(t)}{z\,\tilde{z}\left(\left(1 - t \, S_W^0(t)\right)^2 - z_t\,z_t'\,\tilde{z}_t\, \tilde{z}_t'\, S_V^0(t)\, S_U^0(t) \, t^2\right)} \\
    \\
    S_W(z, \tilde{z}, t) = \frac{S_W^0(t)\, \left(1 - t \, S_W^0(t)\right) + z_t\,z_t'\,\tilde{z}_t\, \tilde{z}_t'\, S_V^0(t)\, S_U^0(t) \, t}{\left(1 - t \, S_W^0(t)\right)^2 - z_t\,z_t'\,\tilde{z}_t\, \tilde{z}_t'\, S_V^0(t)\, S_U^0(t) \, t^2} \,,
\end{cases}
$$
where $z_t\,$, $z_t'\,$, $\tilde{z}_t$ and $\tilde{z}_t'$ are defined in (\ref{eq:G_implicit}) and (\ref{eq:G_tilde_implicit}). For the detailed resolution, see Appendix \ref{subsec:solving_system}. We stress that identifying this limiting differential system and solving it is the most crucial part of our work.

Since $S_V^{(N)}$ converges to a deterministic limit $S_V\,$, we deduce that it is also the limit of its mean. The eigenvalues being deterministic in the scaling limit (we expect them to stick to the quantiles of their limiting deterministic distribution), the expectation is asymptotically taken only on the overlaps, meaning we have  
$$
S_V(z,\tilde{z}, t) = \alpha \, \int_{\mathbb{R}} \, \int_{\mathbb{R}} \, \frac{\bar{V}(\mu, \lambda, t)\, \tilde{\rho}(\mu,t) \, \rho(\lambda,t)}{(\tilde{z} - \mu)(z - \lambda)}\, d\mu \, d\lambda \,.
$$
Therefore, we can recover $\bar{V}$ from $S_V$ using the inversion formula derived in \cite{bun2018overlaps} and used in \cite{attal2024interlacing}\,,
$$
\bar{V}(\mu, \lambda, t) = \lim_{\varepsilon \to 0^+}\frac{1}{2 \pi^2 \alpha \, \tilde{\rho}(\mu,t)\, \rho(\lambda,t)}\, \Re\left[S_V(\lambda - i \, \varepsilon, \mu + i\, \varepsilon,  t) - S_V(\lambda - i \, \varepsilon, \mu - i \, \varepsilon, t)\right] \,,
$$
for any $\mu$ in the support of $\tilde{\rho}(\cdot,t)$ and $\lambda$ in the support of $\rho(\cdot,t)\,$.
Similarly, we have 
$$
S_W(z, \tilde{z}, t) = \alpha \, \int_{\mathbb{R}}\, \int_{\mathbb{R}}\, \frac{\sqrt{\mu \lambda}\, \bar{W}(\mu, \lambda, t)\, \tilde{\rho}(\mu, t)\, \rho(\lambda , t)}{(\tilde{z} - \mu)(z - \lambda)}\, d\mu \, d\lambda\,,
$$
and 
$$
\bar{W}(\mu, \lambda , t) = \lim_{\varepsilon \to 0^+}\frac{1}{2 \pi^2 \alpha \, \sqrt{\mu \lambda}\, \tilde{\rho}(\mu,t)\, \rho(\lambda,t)}\, \Re\left[S_W(\lambda - i \, \varepsilon, \mu + i\, \varepsilon,  t) - S_W(\lambda - i \, \varepsilon, \mu - i \, \varepsilon, t)\right]\,.
$$
The case of $\bar{U}$ is a bit trickier, as we need to split $S_U$ into four parts. Indeed, one has in the scaling limit 
\begin{align*}
    S_U(z, \tilde{z}, t) &= \alpha \,\int_{\mathbb{R}}\, \int_{\mathbb{R}}\, \frac{\bar{U}(\mu, \lambda, t)\, \tilde{\rho}(\mu,t)\, \rho(\lambda, t)}{(\tilde{z} - \mu)(z - \lambda)}\, d\mu \, d\lambda + \frac{\frac{\beta}{q} - \alpha}{\tilde{z}}\, \int_{\mathbb{R}}\, \frac{\bar{U}(0, \lambda, t)\, \rho(\lambda, t)}{(z - \lambda)}\, d\lambda \\
    &\quad + \alpha\, \frac{\frac{1}{q} - 1}{z}\, \int_{\mathbb{R}}\, \frac{\bar{U}(\mu, 0, t)\, \tilde{\rho}(\mu, t)}{(\tilde{z} - \mu)}\, d\mu + \frac{(\frac{\beta}{q} - \alpha)\, (\frac{1}{q} - 1)}{z \tilde{z}}\, \bar{U}(0, 0, t) \,.
\end{align*}
Therefore, we need to use four different inversion formulas to extract $\bar{U}$ in each case:
\begin{itemize}
    \item[$\bullet$] For $\mu$ in the support of $\tilde{\rho}(\cdot, t)$ and $\lambda$ in the support of $\rho(\cdot, t)\,$, we use the same inversion than for $S_V$ and $S_W\,$,
    $$
\bar{U}(\mu, \lambda , t) = \lim_{\varepsilon \to 0^+}\frac{1}{2 \pi^2 \alpha \, \tilde{\rho}(\mu,t)\, \rho(\lambda,t)}\, \Re\left[S_U(\lambda - i \, \varepsilon, \mu + i\, \varepsilon,  t) - S_U(\lambda - i \, \varepsilon, \mu - i \, \varepsilon, t)\right]\,.
    $$
    \item[$\bullet$] For $\mu = 0$ and $\lambda$ in the support of $\rho(\cdot,t)\,$, we use the classical Sokhotski-Plemelj formula already introduced for the Stieltjes transforms (\ref{eq:Stieltjes_inversion_formula}),
    $$
    \bar{U}(0, \lambda, t) = \lim_{\varepsilon \to 0^+}\, \frac{1}{\pi\, (\frac{\beta}{q} - \alpha) \,\rho(\lambda ,t)}\, \Im\left[i\,\varepsilon \, S_U(\lambda - i \,\varepsilon, i \,\varepsilon, t)\right]\,.
    $$
    Taking $\tilde{z} = i\,\varepsilon\,$, multipliying it with $S_U\,$, and sending $\varepsilon$ to $0$ causes all other integrals to vanish because the supports of $\rho$ and $\tilde{\rho}$ are included in $\mathbb{R}_+^*\,$.
    \item[$\bullet$] We use the same method for $\mu$ in the support of $\tilde{\rho}(\cdot, t)$ and $\lambda = 0\,$,
    $$
    \bar{U}(\mu, 0, t) = \lim_{\varepsilon \to 0^+} \, \frac{1}{\pi \, \alpha \, (\frac{1}{q} - 1)\, \tilde{\rho}(\mu, t)}\, \Im \left[i \, \varepsilon \, S_U(i \, \varepsilon, \mu - i \, \varepsilon, t) \right]\,.
    $$
    \item[$\bullet$] The last case is simpler, we take $z = \tilde{z} = i \, \varepsilon$ and send $\varepsilon$ to $0$ which gives
    $$
    \bar{U}(0, 0, t) = \lim_{\varepsilon \to 0^+} \frac{1}{(\frac{\beta}{q} - \alpha)\, (\frac{1}{q} - 1)}\, (i\, \varepsilon)^2\, S_U(i\, \varepsilon, i \, \varepsilon, t)\,.
    $$
\end{itemize}

We are now ready to state our formulas for $\bar{V}\,,\bar{U}$ and $\bar{W}\,$ for a general initial condition $A\,$, from which one can always compute $S_V(\cdot, \cdot, 0)\,, S_U(\cdot, \cdot, 0)$ and $S_W(\cdot, \cdot, 0)\,$. Using the notations $y_t := 1 - t \, v(\lambda, t) - i \pi t \, \rho(\lambda, t)\,$, $y_t' := \lambda \, y_t - \left(q^{-1} - 1\right)\, t\,$, $\tilde{y}_t := 1 - \alpha t \, \tilde{v}(\mu,t) - i \alpha \pi t \, \tilde{\rho}(\mu,t)\,$ and $\tilde{y}_t' := \mu \, \tilde{y}_t - \left(\beta \,/\, q - \alpha\right)\, t\,,$ along with $x_A(t) := S_x\left(y_t \, y_t' \,, \tilde{y}_t \, \tilde{y}_t', 0\right)$ and $x_A^*(t) := S_x\left(y_t \, y_t' \,, \tilde{y}_t^* \, (\tilde{y}_t')^*, 0\right)\,$ for $x \in \{V\,,U\,,W\}\,$, we have the following explicit formulas for $\mu\,,\lambda > 0$ (more precisely in the respective supports of $\tilde{\rho}(\cdot, t)$ and $\rho(\cdot, t)$):
\vspace{0.5cm}
\begin{equation}
    \label{eq:General_solution}
\begin{cases}
    \bar{V}(\mu, \lambda , t) = \frac{1}{Z}\,\Re\left[\frac{y_t \,\tilde{y}_t^*\, V_{A}^*}{\left(1 - t \,W_A^*\right)^2 - y_t\,y_t'\, \tilde{y}_t^* \, (\tilde{y}_t')^* \, V_A^* \, U_A^* \, t^2} -  \frac{y_t\, \tilde{y}_t \,V_A}{\left(1 - t\, W_A\right)^2 - y_t\,y_t'\, \tilde{y}_t \, \tilde{y}_t' \, V_A \, U_A \, t^2} \right] \,,\\
    \\
    \bar{U}(\mu, \lambda ,t ) = \frac{1}{\mu \, \lambda \, Z} \, \Re\left[\frac{y_t' \,(\tilde{y}_t')^*\, U_A^*}{\left(1 - t \,W_A^*\right)^2 - y_t\,y_t'\, \tilde{y}_t^* \, (\tilde{y}_t')^* \, V_A^* \, U_A^* \, t^2} -  \frac{y_t' \,\tilde{y}_t' \,U_A}{\left(1 - t\, W_A\right)^2 - y_t\,y_t'\, \tilde{y}_t \, \tilde{y}_t' \, V_A \, U_A \, t^2}\right]\,,\\
    \\
    \bar{W}(\mu, \lambda , t) = \frac{1}{\sqrt{\mu \lambda}\, Z}\, \Re\left[\frac{W_A^* \, \left(1 - t \, W_A^*\right) + y_t\, y_t' \, \tilde{y}_t^* \, (\tilde{y}_t')^* \, V_A^*\,U_A^*\, t}{\left(1 - t \,W_A^*\right)^2 - y_t\,y_t'\, \tilde{y}_t^* \, (\tilde{y}_t')^* \, V_A^* \, U_A^* \, t^2} -  \frac{W_A\, \left(1 - t \, W_A\right) + y_t \, y_t' \, \tilde{y}_t \, \tilde{y}_t'\, V_A \, U_A \, t}{\left(1 - t\, W_A\right)^2 - y_t\,y_t'\, \tilde{y}_t \, \tilde{y}_t' \, V_A \, U_A \, t^2}\right]\,.
\end{cases}
\vspace{0.5cm}
\end{equation}
where $Z$ is the normalisation $2 \alpha \pi^2 \, \tilde{\rho}(\mu,t)\, \rho(\lambda,t)\,$. For the other cases, we introduce $G_t = lim_{\varepsilon \to 0^+}\, G(i\,\varepsilon, t) = - \int_{\mathbb{R}}\, \frac{\rho(\lambda, t)}{\lambda}\, d\lambda$ as well as $\tilde{G}_t = - \int_{\mathbb{R}}\, \frac{\tilde{\rho}(\mu, t)}{\mu}\, d\mu\,$. Setting $c := \frac{1}{q} - 1$ and $\tilde{c} := \frac{\beta}{q} - \alpha\,$ we have

\begin{align*}
    &\bar{U}(0,\lambda, t) = \frac{1}{\pi \, \lambda\, \rho(\lambda , t)} \, \Im\left[\frac{-y_t'\, U_A(t)\,t}{(1 - t\,W_A(t))^2 + \tilde{c}\, y_t\,y_t'\,(1 - \alpha t\, \tilde{G}_t)\,V_A(t)\,U_A(t)\,t^3}\right]\,,
\end{align*}
where $x_A(t) := S_x\left(y_t\, y_t'\,, - (1 - \alpha t\,\tilde{G}_t)\,\tilde{c}\,t\,, 0\right)$ for $x \in \{V\,, U\,, W\}\,$,
\vspace{0.5cm}
$$
\bar{U}(\mu, 0, t) = \frac{1}{\pi \, \alpha\, \mu\, \tilde{\rho}(\mu, t)}\, \Im\left[\frac{- \tilde{y}_t'\,U_A(t)\,t}{(1 - t\, W_A(t))^2 + c \, \tilde{y}_t\, \tilde{y}_t'\, (1 - t \, G_t)\, V_A(t)\, U_A(t)\, t^3}\right]\,,
$$
where $x_A(t) := S_x\left(- (1 - t\, G_t)\, c\, t\,,\tilde{y}_t\, \tilde{y}_t'\,, 0 \right)$ for $x \in \{V\,, U\,, W\}\,$,
\vspace{0.5cm}
$$
\bar{U}(0,0,t) = \frac{U_A(t)\, t^2}{(1 - t\, W_A(t))^2 - c \, \tilde{c}\, V_A(t)\, U_A(t)\, t^4}\,,
$$
where $x_A(t) := S_x\left(- (1 - t\, G_t)\, c\,t\,, - (1 - \alpha t \, \tilde{G}_t)\, \tilde{c}\, t\,, 0\right)$ for $x \in \{V\,, U\,, W\}\,$.

Hence, we have been able to compute the exact limits of $\mathbb{E}\left[N \, \bk{\tilde{v}_i^t}{v_j^t}^2\right]\,$, $\mathbb{E}\left[N \, \bk{\tilde{u}_i^t}{u_j^t}^2\right]$ and $\mathbb{E}\left[N \,\bk{\tilde{u}_i^t}{u_j^t}\, \bk{\tilde{v}_i^t}{v_j^t}\right]$ for eigenvectors in the bulk. These formulas are completely explicit given the initial condition $A\,$.

\subsection{The Marchenko-Pastur Case}
In this subsection, we show that our formulas simplify when $A \equiv 0\,$.  We have already seen in this case that the distributions $\rho$ and $\tilde{\rho}$ have explicit forms as they are Marchenko-Pastur densities. We also know their Hilbert transforms (see Appendix \ref{subsec:inversion_MP}). Furthermore, since all the eigenvalues are null at $t=0\,$, we have 
$$
S_V(z,\tilde{z}, 0) = \frac{\alpha}{z\, \tilde{z}}\,,
$$
$$
S_U(z,\tilde{z}, 0) = \frac{\beta}{q \, z \, \tilde{z}} \,,
$$
$$
S_W(z,\tilde{z}, 0) = 0 \,.
$$
Therefore, we obtain
$$
\begin{cases}
    S_V(z,\tilde{z}, t) = \frac{\alpha\, z_t\, \tilde{z}_t}{z_t\, z_t'\,\tilde{z}_t\, \tilde{z}_t' - \frac{\alpha\beta}{q}\, t^2} \\
    \\
    S_U(z,\tilde{z},t) = \frac{\beta}{q\,z\,\tilde{z}}\, \frac{z_t'\, \tilde{z}_t'}{z_t\, z_t'\,\tilde{z}_t\, \tilde{z}_t' - \frac{\alpha\beta}{q}\, t^2}\\
    \\
    S_W(z,\tilde{z},t) = \frac{\alpha \beta}{q}\, \frac{t}{z_t\, z_t'\,\tilde{z}_t\, \tilde{z}_t' - \frac{\alpha\beta}{q}\, t^2}\,.
\end{cases}
$$
These forms are explicit and we are able to apply the previous inversion formulas to them, to obtain simplified forms for our goal functions $\bar{V}\,$, $\bar{U}$ and $\bar{W}\,$. 

We get for $\mu \in \left[\left(1 - \sqrt{\frac{\beta}{\alpha \, q}}\right)^2\, \alpha \, t \,, \left(1 + \sqrt{\frac{\beta}{\alpha \, q}}\right)^2\, \alpha \, t\right]$ and $\lambda \in \left[\left(1 - \frac{1}{\sqrt{q}}\right)^2\, t \,, \left(1 + \frac{1}{\sqrt{q}}\right)^2 \, t\right]\,$,
\begin{equation}
    \label{eq:result_A_null}
\begin{cases}
\bar{V}(\mu, \lambda, t) = q\,\frac{(1-\alpha)\, t \,\bar{\mu} + \alpha \, (1 - \beta)\, t \, \bar{\lambda} + (1 - \alpha \beta)\, (\alpha + \frac{1}{q})\, t^2}{(1 - \alpha \beta)^2\, t^2 + q\,(\bar{\lambda} - \bar{\mu})\, (\alpha \beta \, \bar{\lambda} - \bar{\mu})}\\

\vspace{0.1cm}\\

\bar{U}(\mu, \lambda, t) = q\,\frac{(1-\beta)\, t \,\bar{\mu} + \beta \, (1 - \alpha)\, t \, \bar{\lambda} + (1 - \alpha \beta)\, (1 + \frac{\beta}{q})\, t^2}{(1 - \alpha \beta)^2\, t^2 + q\,(\bar{\lambda} - \bar{\mu})\, (\alpha \beta \, \bar{\lambda} - \bar{\mu})}\\

\vspace{0.1cm}\\

\bar{W}(\mu, \lambda, t) = q\,\frac{(1 - \alpha \beta)\, t \, \sqrt{\lambda\, \mu}}{(1 - \alpha \beta)^2\, t^2 + q\,(\bar{\lambda} - \bar{\mu})\, (\alpha \beta \, \bar{\lambda} - \bar{\mu})} \,,
\end{cases}
\end{equation}
where $\bar{\lambda} := \lambda - \left(1 + \frac{1}{q}\right)\, t$ and $\bar{\mu} := \mu - \left(\alpha + \frac{\beta}{q}\right)\, t\,$. This is the most important result of our paper. The calculations leading to these simplifications can be found in Appendix \ref{subsec:inversion_MP}. We note that these are Cauchy-like functions in $\lambda$ or $\mu\,$, as observed in the Wigner setup of \cite{attal2024interlacing}, as well as in \cite{allez2014eigenvectors} and \cite{pacco2023overlaps}. We made our computations in the case $M \geq N$ and $m \geq n\,$, but these three expressions are still valid in any other case. Moreover, they are not affected by a specific choice of bases for the null spaces as they correspond to limiting overlaps between singular vectors associated with non-zero singular values.

\begin{figure}[h]
    \label{fig:Overlaps}
    \centering
    % First image
    \begin{minipage}{0.3\textwidth}
        \centering
        \includegraphics[width=1.\textwidth]{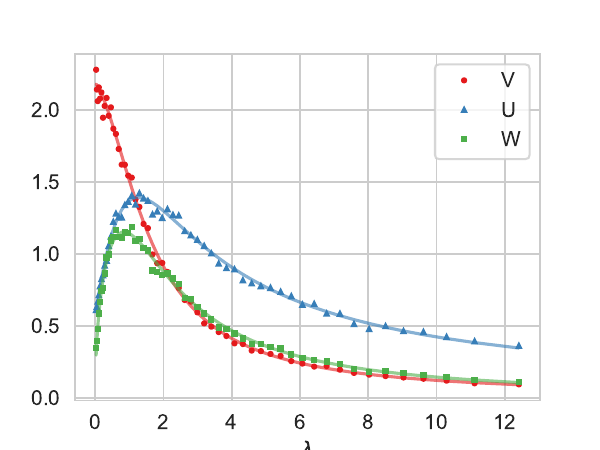}
    \end{minipage}
    \hspace{0.1mm}
    % Second image
    \begin{minipage}{0.3\textwidth}
        \centering
        \includegraphics[width=1.\textwidth]{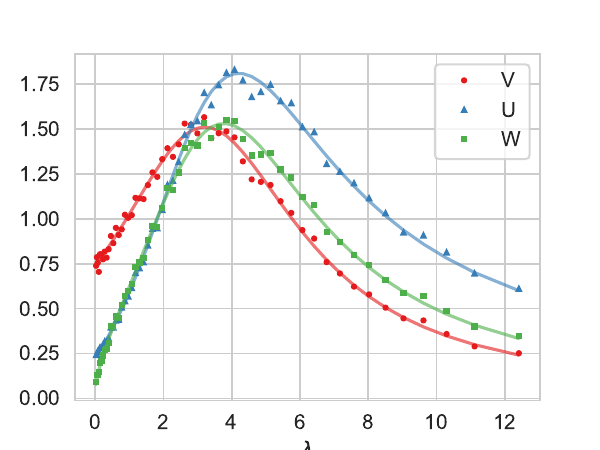}
    \end{minipage}
    \hspace{0.1mm}
    % Third image
    \begin{minipage}{0.3\textwidth}
        \centering
        \includegraphics[width=1.\textwidth]{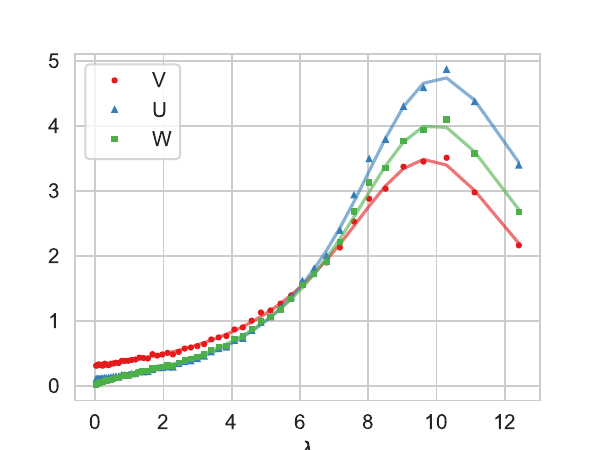}
    \end{minipage}
    \caption{Comparison of our formulas for $\bar{V}\,,\bar{U}$ and $\bar{W}$ with numerical simulations of $N \, \mathbb{E}\left[V_{ij}(t)\right]$ (red plain curve for theory and red circles for data), $N \, \mathbb{E}\left[U_{ij}(t)\right]$ (blue plain curve for theory and blue triangles for data) and $N \, \mathbb{E}\left[W_{ij}(t)\right]$ (green plain curve for theory and green squares for data) for $M = 300\,, q = 0.9\,, \alpha = 0.4\,, \beta = 0.8$ and $t = 3\,$ as a function of $\lambda$ for a fixed $\mu = \mu(x,t)\,$. \textbf{Left:} $x = 0.9\,$. \textbf{Middle:} $x = 0.5\,$. \textbf{Right:} $x = 0.1\,$.}
\end{figure}

Figure \ref{fig:Overlaps} shows a comparison of these formulas with simulated rescaled mean squared overlaps. The fit is excellent.

The other cases for $\bar{U}$ are also simplified into 
\begin{equation}
    \label{eq:result_A_null_ker}
    \begin{cases}
        \bar{U}(0, \lambda , t) = \frac{(1 - \alpha)\, t}{\alpha \, \lambda + (1 - \alpha)\, (\frac{1}{q} - \alpha)\, t}\\
        \\
        \bar{U}(\mu, 0, t) = \frac{(1 - \beta)\, t}{\mu + (1 - \beta)\, (\frac{1}{q} - \alpha)\, t}\\
        \\
        \bar{U}(0,0, t) = \frac{q}{1 - \alpha \, q} \,,
    \end{cases}
\end{equation}
using $\lim_{\varepsilon \to 0^+}\, G(i\, \varepsilon, t) = - \frac{q}{(1 - q)\, t}$ and $\lim_{\varepsilon \to 0^+}\, \tilde{G}(i\, \varepsilon, t) = - \frac{q}{(\beta - \alpha \, q)\, t}$ in the Marchenko-Pastur setup. The forms (\ref{eq:result_A_null_ker}) are specific to our choice of structure for the null spaces made in the introduction, and to the situation $M \geq N\,, m \geq n\,$. Note that numerically there can be some differences with the overlaps obtained with simulation for certain choices of the parameters $q\,, \alpha$ and $\beta\,$, due to the finite matrix size and unexchangeability of the singular vectors.

\newpage
\appendix
\section*{Appendices}
\renewcommand{\thesubsection}{\Alph{subsection}}
\renewcommand{\theequation}{\Alph{subsection}.\arabic{equation}}
\numberwithin{equation}{subsection}

\subsection{Correlation}
\label{subsec:App_correlation}

Let $1 \leq i \leq m \,$, $1 \leq l \leq n\,$, $1 \leq j \leq M\,$ and $1 \leq k \leq N\,$, we have:
\begin{align*}
\bk{u_j^t}{dB_t \, v_k^t} \, \bk{\tilde{u}_i^t}{d\tilde{B}_t \, \tilde{v}_l^t} &= \left( \sum_{r = 1}^M \, \sum_{s = 1}^N \, u_{jr}^t \, v_{ks}^t \, dB_t^{rs} \right) \left( \sum_{r=1}^m \, \sum_{s=1}^n \, \tilde{u}_{ir}^t \, \tilde{v}_{ls}^t dB_t^{rs} \right) \\
&= \sum_{r=1}^m \, \sum_{s=1}^n \, u_{jr}^t \, \tilde{u}_{ir}^t \, v_{ks}^t \, \tilde{v}_{ls}^t \, dt \,.
\end{align*}
We recall that for $1 \leq l \leq n\,$, $\tilde{v}_{ls}^t = 0$ if $s > n$ and for $1\leq i \leq m\,$, $\tilde{u}_{ir}^t = 0$ if $r > m\,$. Thus we indeed have:
\begin{equation*}
\bk{u_j^t}{dB_t \, v_k^t} \, \bk{\tilde{u}_i^t}{d\tilde{B}_t \, \tilde{v}_l^t} = \bk{\tilde{u}_i^t}{u_j^t} \, \bk{\tilde{v}_l^t}{v_k^t} \, dt \,.
\end{equation*}

\subsection{Burgers Equation}
\label{subsec:App_burgers}
\subsubsection{Deriving the Equation}
\label{subsubsec:App_deriving_burgers}
Applying Itô's lemma gives:
\begin{align*}
dG_N(z,t) &= \frac{2}{N \sqrt{N}} \, \sum_{j = 1}^N \, \frac{\sqrt{\lambda_j^t} \, db_j(t)}{(z - \lambda_j^t)^2} + \frac{M}{N^2} \sum_{\substack{j = 1}}^N \, \frac{1}{(z - \lambda_j^t)^2} \, dt \\
&\quad +  \frac{1}{N^2} \, \sum_{\substack{j,k = 1 \\ k \neq j}}^N \, \frac{\lambda_j^t + \lambda_k^t}{(\lambda_j^t - \lambda_k^t) (z-\lambda_j^t)^2} \, dt + \frac{4}{N^2} \, \sum_{j=1}^N\, \frac{\lambda_j^t}{(z - \lambda_j^t)^3} \, dt\,.
\end{align*}
The first and last sum go to 0 in the scaling limit, and the second one converges to $ - \frac{\partial_z G(z,t)}{q} \, dt\,$. We need to perform some manipulations to deal with the third sum, that we denote by $\Sigma\,dt\,$. We first split it into $\Sigma\,/\,2 + \Sigma \,/\,2$ and invert the indices in the second term. Regrouping the two sums and applying the identity
\begin{equation*}
    \frac{1}{(z - b_j)^2} - \frac{1}{(z - b_i)^2} = \frac{(b_j - b_i)(2z - b_j - b_i)}{(z - b_j)^2 (z-b_i)^2} \,,
\end{equation*} 
we get 
\begin{align*}
\Sigma &= \frac{1}{2N^2}\, \sum_{\substack{j,k = 1 \\ k \neq j}}^{N}\, \frac{(\lambda_j^t + \lambda_k^t)\, (2z - \lambda_j^t - \lambda_k^t)}{(z - \lambda_j^t)^2 (z - \lambda_k^t)^2}\\
&= \frac{1}{2N^2}\, \sum_{\substack{j,k = 1 \\ k \neq j}}^{N}\, \frac{\lambda_j^t + \lambda_k^t}{(z - \lambda_j^t)(z - \lambda_k^t)^2} + \frac{1}{2N^2}\, \sum_{\substack{j,k = 1 \\ k \neq j}}^{N}\, \frac{\lambda_j^t + \lambda_k^t}{(z - \lambda_j^t)^2(z - \lambda_k^t)} \\
&= \frac{1}{N^2}\, \sum_{\substack{j,k = 1 \\ k \neq j}}^{N}\, \frac{\lambda_j^t + \lambda_k^t}{(z - \lambda_j^t)(z - \lambda_k^t)^2} \,.
\end{align*}
We split this forms into two sums:
$$
\Sigma = \frac{1}{N^2}\, \sum_{\substack{j,k = 1 \\ k \neq j}}^{N}\, \frac{\lambda_j^t}{(z - \lambda_j^t)(z - \lambda_k^t)^2} + \frac{1}{N^2}\, \sum_{\substack{j,k = 1 \\ k \neq j}}^{N}\, \frac{\lambda_k^t}{(z - \lambda_j^t)(z - \lambda_k^t)^2}\,,
$$
where the first sum equals, using $\lambda \, / \, (z - \lambda) = -1 + z \,/\, (z-\lambda)\,$ and adding the missing diagonal terms,
\begin{align*}
\left( -1 + \frac{z}{N} \, \sum_{j = 1}^N \, \frac{1}{z - \lambda_j^t} \right)\left(\frac{1}{N} \, \sum_{k=1}^N \, \frac{1}{(z - \lambda_k^t)^2}\right) - \frac{1}{N^2} \, \sum_{j= 1}^N \, \frac{\lambda_j^t}{(z - \lambda_j^t)^3} \,,
\end{align*}
which converges to
$$
\left(1 - z \, G(z,t) \right) \, \partial_z G(z,t) \,,
$$
and the second sum equals
$$
\left(\frac{1}{N} \, \sum_{j=1}^N \, \frac{1}{z - \lambda_j^t}\right) \left(- \frac{1}{N} \, \sum_{k = 1}^N \, \frac{1}{z - \lambda_k^t} + \frac{z}{N} \, \sum_{k=1}^N \, \frac{1}{(z - \lambda_k^t)^2} \right) - \frac{1}{N^2} \, \sum_{j= 1}^N \, \frac{\lambda_j^t}{(z - \lambda_j^t)^3} \,,
$$
which converges to
$$
G(z,t) \, \left(-G(z,t) - z\, \partial_zG(z,t) \right) \,.
$$
Finally, regrouping all the terms leads to the announced limiting equation (\ref{eq:G_Burgers})
$$
\partial_t G(z,t) = \left(1 - \frac{1}{q} - 2z \, G(z,t) \right) \, \partial_zG(z,t) - G^2(z,t) \,.
$$

\subsubsection{Solving the Equation}
\label{subsubsec:App_solving_burgers}
In order to obtain the implicit equation (\ref{eq:G_implicit}) satisfied by $G\,$, we use the method of characteristics. We introduce two functions of a new variable $s\,$: $z(s)$ and $t(s)\,$. We define $\hat{G}(s) := G(z(s),t(s))\,$, so that the chain rule gives us
\begin{align*}
\frac{d\hat{G}}{ds} &= \partial_z G(z(s),t(s)) \, \frac{dz}{ds} + \partial_t G(z(s),t(s)) \, \frac{dt}{ds} \\
&= \left(\frac{dz}{ds} + \left(1 - \frac{1}{q} - 2z(s) \, \hat{G}\right) \, \frac{dt}{ds} \right) \, \partial_z G(z(s), t(s)) - \hat{G}^2 \, \frac{dt}{ds} \,.
\end{align*}
Therefore, if we choose the functions $z$ and $t$ such that
\begin{equation*}
\begin{cases}
\frac{dt}{ds} = 1 \\
\frac{dz}{ds} = 2z(s) \, \hat{G}(s) + \frac{1}{q} - 1 \,,
\end{cases}
\end{equation*}
then $d\hat{G} \,/\, ds = - \hat{G}^2\,$, meaning $\hat{G}(s) = \hat{G}(0) \,/\, (1 + s \,\hat{G}(0))\,$. This simplifies the differential equation on $z$ which allows us to obtain
$$
z(s) = \left(1 + s \,\hat{G}(0)\right) \, \left(z(0) \, \left(1 + s \,\hat{G}(0)\right) + \left(\frac{1}{q} - 1\right)\,s\right)\,.
$$
Finally, the solution $\hat{G}$ gives $\hat{G}(0) = \hat{G}(s) \,/\, (1 - s \,\hat{G}(s))\,$, i.e. for any $s\,$,
$$
G(z(0), t(0)) = \frac{G(z(s), t(0) + s)}{1 - s \,G(z(s), t(0) + s)}\,.
$$
Evaluating this at $s = -t(0)$ and noticing $z(0)$ and $t(0)$ are free parameters, we obtain the announced implicit equation (\ref{eq:G_implicit})
$$
G(z,t) = \frac{G\left((1 - t\, G(z,t))\,\left(z\,(1-t\,G(z,t)) + \left(1 - q^{-1}\right)\,t\right)\,, 0 \right)}{1 + t \, G\left((1 - t\, G(z,t))\,\left(z\,(1-t\,G(z,t)) + \left(1 - q^{-1}\right)\,t\right)\,, 0 \right)} \,.
$$

\subsection{Itô Dynamics of the Squared Overlaps}
\label{subsec:App_ito_overlaps}
We compute here the Itô dynamics of the different squared overlaps $V_{ij}(t)\,$, $U_{ij}(t)$ and $W_{ij}(t)\,$. We detail the calculations for $V_{ij}$ and state the dynamics for the other cases, as the calculations are similar. For readability, we use the notation $[\,\cdot\,]$ by
\begin{equation}
\label{eq:def_crochet}
\left[a_{iljk}\right] := a_{iljk} + a_{ilkj} + a_{lijk} + a_{likj} \,.
\end{equation}
In addition, we define $\ov{i}{l}{j}{k} := \bk{\tilde{v}_i^t}{v_j^t} \, \bk{\tilde{v}_l^t}{v_k^t}\,, \ou{i}{l}{j}{k} := \bk{\tilde{u}_i^t}{u_j^t}\, \bk{\tilde{u}_l^t}{u_k^t}$ and $\ow{i}{l}{j}{k} := \bk{\tilde{v}_i^t}{v_j^t}\, \bk{\tilde{u}_l^t}{u_k^t}\,$.

Let $1 \leq i \leq n$ and $1 \leq j \leq N\,$, we first compute
\begin{align*}
d\bk{\tilde{v}_i^t}{v_j^t} &= \bk{d\tilde{v}_i^t}{v_j^t} + \bk{\tilde{v}_i^t}{dv_j^t} + \bk{d\tilde{v}_i^t}{dv_j^t} \\
&= -\frac{1}{2N} \, \sum_{\substack{l = 1 \\ l  \neq i}}^n \, \frac{\mu_i^t + \mu_l^t}{(\mu_i^t - \mu_l^t)^2} \, \bk{\tilde{v}_i^t}{v_j^t} \, dt + \frac{1}{\sqrt{N}} \, \sum_{\substack{l =1 \\ l \neq i}}^n \, \frac{\sqrt{\mu_i^t} \, d\tilde{w}_{il}(t) + \sqrt{\mu_l^t}\, d\tilde{w}_{li}(t)}{\mu_i^t - \mu_l^t}\, \bk{\tilde{v}_l^t}{v_j^t} \\
&\quad -\frac{1}{2N} \, \sum_{\substack{k = 1 \\ k \neq j}}^N \, \frac{\lambda_j^t + \lambda_k^t}{(\lambda_j^t - \lambda_k^t)^2}\, \bk{\tilde{v}_i^t}{v_j^t}\, dt + \frac{1}{\sqrt{N}} \, \sum_{\substack{k = 1 \\ k \neq j}}^N \, \frac{\sqrt{\lambda_j^t}\, dw_{jk}(t) + \sqrt{\lambda_k^t}\, dw_{kj}(t)}{\lambda_j^t - \lambda_k^t} \, \bk{\tilde{v}_i^t}{v_k^t} \\
&\quad + \frac{1}{N} \, \sum_{\substack{l = 1 \\ l \neq i}}^n \, \sum_{\substack{k = 1 \\ k \neq j}}^N \, \frac{A_{iljk}^t}{(\mu_i^t - \mu_l^t)(\lambda_j^t - \lambda_k^t)} \, \bk{\tilde{v}_l^t}{v_k^t} \,,
\end{align*}
where for any $l \neq i$ in $\{1\,;...\,;n\}$ and any $k \neq j$ in $\{1\,;...\,;N\}\,$,
\begin{align*}
A_{iljk}^t &:= \left(\sqrt{\mu_i^t} \, d\tilde{w}_{il}(t) + \sqrt{\mu_l^t} \, d\tilde{w}_{li}(t)\right) \left( \sqrt{\lambda_j^t} \, dw_{jk}(t) + \sqrt{\lambda_k^t} \, dw_{kj}(t) \right) \\
&= \left[\sqrt{\mu_i^t \lambda_j^t} \, \ow{l}{i}{k}{j}\right] \, dt \,.
\end{align*}
Now, we can compute the dynamics of the squared overlaps:
\begin{align*}
dV_{ij}(t) &= 2 \bk{\tilde{v}_i^t}{v_j^t} \, d\bk{\tilde{v}_i^t}{v_j^t} + \left(d\bk{\tilde{v}_i^t}{v_j^t}\right)^2 \\
&= -\frac{1}{N} \, \sum_{\substack{l = 1 \\ l \neq i}}^n \, \frac{\mu_i^t + \mu_l^t}{(\mu_i^t- \mu_l^t)^2}\, V_{ij} \, dt + \frac{2}{\sqrt{N}}\, \sum_{\substack{l = 1 \\ l \neq i}}^n \, \frac{\sqrt{\mu_i^t}\, d\tilde{w}_{il}(t) + \sqrt{\mu_l^t}\, d\tilde{w}_{li}(t)}{\mu_i^t - \mu_l^t} \, \bk{\tilde{v}_l^t}{v_j^t} \, \bk{\tilde{v}_i^t}{v_j^t} \\
&\quad - \frac{1}{N} \, \sum_{\substack{k = 1 \\ k \neq j}}^N \, \frac{\lambda_j^t + \lambda_k^t}{(\lambda_j^t - \lambda_k^t)^2} \, V_{ij}\, dt + \frac{2}{\sqrt{N}} \, \sum_{\substack{k = 1 \\ k \neq j}}^N \, \frac{\sqrt{\lambda_j^t} \, dw_{jk}(t) + \sqrt{\lambda_k^t} \, dw_{kj}(t)}{\lambda_j^t - \lambda_k^t} \, \bk{\tilde{v}_i^t}{v_k^t} \, \bk{\tilde{v}_i^t}{v_j^t} \\
&\quad + \frac{2}{N} \, \sum_{\substack{l = 1 \\ l \neq i}}^n \, \sum_{\substack{k = 1 \\ k \neq j}}^N \, \frac{A_{iljk}^t}{(\mu_i^t - \mu_l^t)(\lambda_j^t - \lambda_k^t)} \, \bk{\tilde{v}_l^t}{v_k^t} \, \bk{\tilde{v}_i^t}{v_j^t}
\end{align*}
\begin{align*}
&\hspace{-3.5cm}\quad + \frac{1}{N} \, \sum_{\substack{l=1 \\ l \neq i}}^n \, \frac{\mu_i^t + \mu_l^t}{(\mu_i^t - \mu_l^t)^2} \, V_{lj}\,dt + \frac{1}{N} \, \sum_{\substack{k = 1 \\ k \neq j}}^N \, \frac{\lambda_j^t + \lambda_k^t}{(\lambda_j^t - \lambda_k^t)^2} \, V_{ik} \, dt \\
&\hspace{-3.5cm}\quad + \frac{2}{N} \, \sum_{\substack{l = 1 \\ l \neq i}}^n \, \sum_{\substack{k = 1 \\ k \neq j}}^N \, \frac{A_{iljk}^t}{(\mu_i^t - \mu_l^t)(\lambda_j^t - \lambda_k^t)} \, \bk{\tilde{v}_l^t}{v_j^t} \, \bk{\tilde{v}_i^t}{v_k^t} \,,
\end{align*}
which can be rewritten as
\begin{align*}
dV_{ij}(t) &= \frac{1}{N} \, \sum_{\substack{l = 1 \\ l \neq i}}^n \, \frac{\mu_i^t + \mu_l^t}{(\mu_i^t - \mu_l^t)^2} \, \left(V_{lj} - V_{ij}\right)\, dt + \frac{1}{N} \, \sum_{\substack{k = 1 \\ k  \neq j}}^N \, \frac{\lambda_j^t + \lambda_k^t}{(\lambda_j^t - \lambda_k^t)^2} \, \left(V_{ik} - V_{ij}\right)\,dt \\
&\quad + \frac{2}{N} \, \sum_{\substack{l = 1 \\ l \neq i}}^n \, \sum_{\substack{k = 1 \\ k \neq j}}^N \, \frac{\left[V_{ij}\,\bar{W}_{lk}\right] + \left[\sqrt{\mu_i^t \lambda_j^t} \, \ov{i}{l}{k}{j} \, \ow{l}{i}{k}{j} \right]}{(\mu_i^t - \mu_l^t)(\lambda_j^t - \lambda_k^t)} \, dt \\
& \quad + \frac{2}{\sqrt{N}} \, \sum_{\substack{l = 1 \\ l \neq i}}^n \, \frac{\sqrt{\mu_i^t}\,d\tilde{w}_{il}(t) + \sqrt{\mu_l^t}\, d\tilde{w}_{li}(t)}{\mu_i^t - \mu_l^t} \, \ov{i}{l}{j}{j} \\
& \quad + \frac{2}{\sqrt{N}} \, \sum_{\substack{k =1\\ k \neq j}}^N \, \frac{\sqrt{\lambda_j^t}\, dw_{jk}(t) + \sqrt{\lambda_k^t}\, dw_{kj}(t)}{\lambda_j^t - \lambda_k^t} \, \ov{i}{i}{j}{k} \,,
\end{align*}
where $\bar{W}_{ij}(t) := \sqrt{\mu_i^t \lambda_j^t} \, W_{ij}(t) \,$.

Similarly, one finds that
\begin{align*}
dU_{ij}(t) &= \frac{1}{N} \, \sum_{\substack{l = 1 \\ l \neq i}}^m \, \frac{\mu_i^t + \mu_l^t}{(\mu_i^t - \mu_l^t)^2}\, \left(U_{lj} - U_{ij}\right)\,dt + \frac{1}{N} \, \sum_{\substack{k=1\\ k \neq j}}^M \, \frac{\lambda_j^t + \lambda_k^t}{(\lambda_j^t - \lambda_k^t)^2} \, \left(U_{ik} - U_{ij}\right)\, dt \\
&\quad + \frac{2}{N} \, \sum_{\substack{l=1 \\ l \neq i}}^m \, \sum_{\substack{k=1\\ k \neq j}}^M\, \frac{\left[U_{ij}\,\bar{W}_{lk}\right] + \left[\sqrt{\mu_i^t \lambda_j^t}\, \ou{i}{l}{k}{j} \, \ow{i}{l}{j}{k}\right]}{(\mu_i^t - \mu_l^t)(\lambda_j^t - \lambda_k^t)}\, dt \\
&\quad + \frac{2}{\sqrt{N}} \, \sum_{\substack{l= 1 \\ l \neq i}}^m\, \frac{\sqrt{\mu_i^t}\, d\tilde{w}_{li}(t) + \sqrt{\mu_l^t}\, d\tilde{w}_{il}(t)}{\mu_i^t - \mu_l^t} \, \ou{i}{l}{j}{j} \\
&\quad + \frac{2}{\sqrt{N}}\, \sum_{\substack{k= 1 \\ k \neq j}}^M \, \frac{\sqrt{\lambda_j^t}\, dw_{kj}(t) + \sqrt{\lambda_k^t}\, dw_{jk}(t)}{\lambda_j^t - \lambda_k^t}\, \ou{i}{i}{j}{k} \,.
\end{align*}

Finally, for $W_{ij}\,$, the form is quite heavy as we mix sums with indices ending at four different bounds: $n\,$, $m\,$, $N$ and $M\,$. We find
\begin{align*}
dW_{ij}(t) &= \frac{1}{N} \, \sum_{\substack{l = 1 \\ l \neq i}}^n \, \frac{2\,\sqrt{\mu_i^t \mu_l^t}\, W_{lj} - \left(\mu_i^t + \mu_l^t\right)\, W_{ij}}{(\mu_i^t - \mu_l^t)^2} \, dt + \frac{n - m}{2N\mu_i^t} \, W_{ij} \, dt \\
& \quad + \frac{1}{N} \, \sum_{\substack{k = 1 \\ k \neq j}}^N \, \frac{2\, \sqrt{\lambda_j^t \lambda_k^t}\, W_{ik} - \left(\lambda_j^t + \lambda_k^t\right)\, W_{ij}}{(\lambda_j^t - \lambda_k^t)^2} \, + \frac{N - M}{2N\lambda_j^t}\, W_{ij} \, dt \\
& \quad + \frac{1}{N} \, \sum_{\substack{l = 1 \\ l \neq i}}^n \, \sum_{\substack{k = 1 \\ k \neq j}}^N \, \frac{\sqrt{\mu_i^t \lambda_j^t}\, \left[V_{ij}\,U_{lk}\right] + \sqrt{\mu_l^t\lambda_k^t}\, \left[W_{ij}\, W_{lk}\right]}{(\mu_i^t - \mu_l^t)(\lambda_j^t - \lambda_k^t)}\,dt \\
& \quad + \frac{1}{N} \, \sum_{\substack{l = 1 \\ l \neq i}}^n \, \sum_{\substack{k = 1 \\ k \neq j}}^N \, \frac{\sqrt{\mu_i^t\lambda_k^t}\, \left[ \ov{l}{l}{j}{k} \, \ou{i}{i}{j}{k}   \right] + \sqrt{\mu_l^t\lambda_j^t}\, \left[ \ov{i}{l}{k}{k} \, \ou{i}{l}{j}{j} \right]}{(\mu_i^t - \mu_l^t)(\lambda_j^t - \lambda_k^t)}\,dt \\
&\quad + \frac{1}{N} \, \sum_{l=n+1}^m \, \sum_{\substack{k = 1 \\ k \neq j}}^N \, \frac{\sqrt{\lambda_j^t}\, \left(V_{ij}\,U_{lk} + V_{ik}\, U_{lj}\right) + 2\,\sqrt{\lambda_k^t}\, \ov{i}{i}{j}{k}\, \ou{l}{l}{j}{k}}{\sqrt{\mu_i^t}\,(\lambda_j^t - \lambda_k^t)}\, dt\\
&\quad + \frac{1}{N} \, \sum_{\substack{l = 1 \\ l \neq i}}^n \, \sum_{k=N+1}^M \, \frac{\sqrt{\mu_i^t}\, \left(V_{ij}\, U_{lk} + V_{lj}\, U_{ik}\right) + 2 \, \sqrt{\mu_l^t} \, \ov{i}{l}{j}{j} \, \ou{i}{i}{j}{k}}{(\mu_i^t - \mu_l^t) \, \sqrt{\lambda_j^t}} \, dt \\
&\quad + \frac{1}{N\sqrt{\mu_i^t \lambda_j^t}} \, \sum_{i=n+1}^m\, \sum_{j=N+1}^M \,  V_{ij}\, U_{lk} \, dt + \frac{1}{\sqrt{N}} \, \sum_{\substack{l = 1 \\ l \neq i}}^n \, \frac{\sqrt{\mu_i^t}\, d\tilde{w}_{il}(t) + \sqrt{\mu_l^t}\, d\tilde{w}_{li}(t)}{\mu_i^t - \mu_l^t} \, \ow{l}{i}{j}{j} \\
& \quad + \frac{1}{\sqrt{N}} \, \sum_{\substack{l= 1\\ l \neq i}}^n \, \frac{\sqrt{\mu_i^t}\, d\tilde{w}_{li}(t) + \sqrt{\mu_l^t}\, d\tilde{w}_{il}(t)}{\mu_i^t - \mu_l^t} \, \ow{i}{l}{j}{j} + \frac{1}{\sqrt{N \mu_i^t}} \, \sum_{l=n+1}^m\, d\tilde{w}_{li}(t) \, \ow{i}{l}{j}{j} \\
&\quad + \frac{1}{\sqrt{N}} \, \sum_{\substack{k = 1 \\ k \neq j}}^N \,\frac{\sqrt{\lambda_j^t}\, dw_{jk}(t) + \sqrt{\lambda_k^t}\, dw_{kj}(t)}{\lambda_j^t - \lambda_k^t} \, \ow{i}{i}{k}{j}\\
&\quad  + \frac{1}{\sqrt{N}} \sum_{\substack{k = 1 \\ k \neq j}}^N \, \frac{\sqrt{\lambda_j^t}\, dw_{kj}(t) + \sqrt{\lambda_k^t}\, dw_{jk}(t)}{\lambda_j^t - \lambda_k^t} \, \ow{i}{i}{j}{k}  + \frac{1}{\sqrt{N} \lambda_j^t} \sum_{k = N+1}^M \, dw_{kj}(t) \, \ow{i}{i}{j}{k} \,.
\end{align*}

\subsection{System of Partial Differential Equations on the Double Stieltjes Transforms}
\label{subsec:deriving_system}
\subsubsection{First Properties}
Here, we define certain tools that will be our main manipulations to derive the system of partial differential equations. We make use of the notations introduced in Appendix \ref{subsec:App_ito_overlaps}. We first introduce four symmetrisation properties on sums:
\begin{equation}
\label{eq:S1}
\tag{S1}
\text{If $a_{kp} = a_{pk}\,$,} \quad \sum_{\substack{k,p\\ p \neq k}} \, \frac{a_{kp}}{(b_k - b_p)(z - b_k)} = \frac{1}{2} \, \sum_{\substack{k,p \\ p \neq k}} \, \frac{a_{kp}}{(z - b_k)(z-b_p)} \,,
\end{equation}
\begin{equation}
\label{eq:S2}
\tag{S2}
\sum_{\substack{k,p \\ p \neq k}} \, a_{kp} + a_{pk} = 2 \, \sum_{\substack{k,p \\ p \neq k}} \, a_{kp} \,,
\end{equation}
\begin{equation}
\label{eq:S3}
\tag{S3}
\sum_{\substack{k,p \\ p \neq k}} \, \frac{a_{kp} + a_{pk}}{(b_k - b_p)(z-b_k)} = \sum_{\substack{k,p \\ p \neq k}} \, \frac{a_{kp}}{(z-b_k)(z-b_p)} \,,
\end{equation}
\begin{equation}
\label{eq:S4}
\tag{S4}
\sum_{\substack{i,l\\ l \neq i}} \, \sum_{\substack{j,k\\k \neq j}} \, \frac{\left[ a_{iljk}\right]}{(\mu_i^t - \mu_l^t)(\tilde{z}- \mu_i^t)(\lambda_j^t - \lambda_k^t)(z-\lambda_j^t)} = \sum_{\substack{i,l \\ l \neq i}}\, \sum_{\substack{j,k \\ k \neq j}} \, \frac{a_{iljk}}{(\tilde{z} - \mu_i^t)(\tilde{z} - \mu_l^t)(z - \lambda_j^t)(z - \lambda_k^t)} \,.
\end{equation}
These properties can be easily proved:
\begin{itemize}
\item[$\bullet$] For (\ref{eq:S1}), we separate the left sum $S$ into $S \,/\,2 + S\,/\, 2$ and invert the indices in the second term. Then, we apply the identity
\begin{equation}
\label{eq:I}
\tag{I}
\frac{1}{(b_k - b_p)(z-b_k)} - \frac{1}{(b_k - b_p)(z - b_p)} = \frac{1}{(z - b_k)(z- b_p)} \,.
\end{equation}
\item[$\bullet$] (\ref{eq:S2}) is easily obtained by expanding into two sums and inverting the indices in the second one.
\item[$\bullet$] (\ref{eq:S3}) is an application of the two previous properties. Indeed, $a_{kp} + a_{pk}$ is symmetric so we can use (\ref{eq:S1}) and obtain
$$
\frac{1}{2} \, \sum_{\substack{k,p \\ p \neq k}} \frac{ a_{kp} + a_{pk}}{(z - b_k)(z-b_p)} \,.
$$
Symmetrisation (\ref{eq:S2}) then gives the final result.
\item[$\bullet$] Symmetrisation (\ref{eq:S4}) is an application of (\ref{eq:S3}) to each double sum separately, i.e. to indices $i$ and $l$ and then to $j$ and $k\,$.
\end{itemize}

Finally, we prove a reduction property that exploits the specific structure of a certain type of sum that we will encounter several times in our computation. It shows that despite the fact that this sum appears to be of order $\mathcal{O}(1)$ given the order of magnitude of the overlaps in the bulk ($1 \,/\,\sqrt{N}$), it is in fact going to zero in the scaling limit at least as $1 \,/\,N\,$:
\begin{equation}
\label{eq:R}
\tag{R}
\frac{1}{N^2\sqrt{N}} \, \sum_{\substack{(i,l) \in \tilde{I} }} \, \sum_{\substack{(j,k) \in I}} \, \frac{c_{ij}\, \bk{\tilde{v}_i^t}{v_k^t}\, \bk{\tilde{v}_l^t}{v_j^t}\, \bk{\tilde{v}_l^t}{v_k^t}}{(\tilde{z} - \mu_i^t)^{p_1}(\tilde{z} - \mu_l^t)^{p_2}(z - \lambda_j^t)^{p_3}(z-\lambda_k^t)^{p_4}} = \mathcal{O}\left(\frac{1}{N}\right) \,,
\end{equation}
for any $p_1\,,p_2\,,p_3\,,p_4 \geq 0\,$ and $z\,,\tilde{z} \in \mathbb{C} \setminus \mathbb{R}\,$ and $c_{ij} = \mathcal{O}(1)\,$. With $I$ and $\tilde{I}$ both in $\{1\,;...\,;N\}^2$ such that the summands are well defined.

\begin{proof}
We introduce the notations
$$
a_{ij} := \frac{c_{ij}}{(\tilde{z} - \mu_i^t)^{p_1}(z - \lambda_j^t)^{p_3}} \,, \quad b_{ij} := \sum_{\substack{l \\ (i,l) \in \tilde{I}}}\, \sum_{\substack{k \\ (j,k) \in I}}\, \frac{\bk{\tilde{v}_i^t}{v_k^t}\, \bk{\tilde{v}_l^t}{v_j^t}\, \bk{\tilde{v}_l^t}{v_k^t}}{(\tilde{z}-\mu_l^t)^{p_2}(z-\lambda_k^t)^{p_4}} \,,
$$
that are considered null if the indices $(i,j)$ do not allow the correct definition of the terms. We have
$$
\Sigma = \frac{1}{N^2\sqrt{N}}\, \sum_{i,j}\, a_{ij} \, b_{ij} \,.
$$
Using the Cauchy-Schwarz inequality we get
$$
|\Sigma|^2 \leq \frac{1}{N^5} \left( \sum_{i,j}\, |a_{ij}|^2 \right)\, \left( \sum_{i,j}\, |b_{ij}|^2 \right)\,.
$$
Let us treat both sums separately. First,
\begin{align*}
\sum_{i,j}\, |a_{ij}|^2 &= \sum_{i,j}\, \frac{|c_{ij}|^2}{|\tilde{z} - \mu_i^t|^{2p_1}|z-\lambda_j^t|^{2p_3}} = \mathcal{O}(N^2)\,,
\end{align*}
and secondly,
\begin{align*}
\sum_{i,j}\, |b_{ij}|^2  &\leq \sum_{i=1}^N\, \sum_{j=1}^N\, |b_{ij}|^2 \\
&\leq \sum_{i=1}^N \, \sum_{j=1}^N \, \sum_{\substack{l,l' \\ (i,l)\in \tilde{I} \\ (i,l') \in \tilde{I}}}\, \sum_{\substack{k,k'\\ (j,k)\in I \\ (j,k') \in I}}\, \frac{\bk{\tilde{v}_i^t}{v_k^t}\, \bk{\tilde{v}_l^t}{v_j^t}\, \bk{\tilde{v}_l^t}{v_k^t}\,\bk{\tilde{v}_i^t}{v_{k'}^t}\, \bk{\tilde{v}_{l'}^t}{v_j^t}\, \bk{\tilde{v}_{l'}^t}{v_{k'}^t}}{(\tilde{z} - \mu_l^t)^{p_2}(\tilde{z}^* - \mu_{l'}^t)^{p_2}(z-\lambda_k^t)^{p_4}(z^* - \lambda_{k'}^t)^{p_4}} \,.
\end{align*}
Since $v_1^t\,,...\,,v_N^t$ is an orthonormal basis of $\mathbb{R}^N\,$, we have
$$
\sum_{j=1}^N \, \bk{\tilde{v}_l^t}{v_j^t}\, \bk{\tilde{v}_{l'}^t}{v_j^t} = \bk{\tilde{v}_l^t}{\tilde{v}_{l'}^t} = \delta_{ll'} \,,
$$
and similarly
$$
\sum_{i=1}^N \, \bk{\tilde{v}_i^t}{v_k^t}\, \bk{\tilde{v}_i^t}{v_{k'}^t} = \delta_{kk'}\,,
$$
so that
\begin{align*}
\sum_{i,j}\, |b_{ij}|^2 &\leq \sum_{l,k} \,  \frac{V_{lk}(t)}{|\tilde{z} - \mu_l^t|^{2p_2}|z-\lambda_k^t|^{2p_4}} = \mathcal{O}(N) \,.
\end{align*}
Therefore, we get
$$
|\Sigma|^2 = \mathcal{O}\left(\frac{1}{N^2}\right) \,,
$$
which means $\Sigma = \mathcal{O}(1\,/\,N)\,$.
\end{proof}
Note that this property is also satisfied if we replace the overlaps $\bk{\tilde{v}}{v}$ by $\bk{\tilde{u}}{u}\,$, working with indices in $\{1\,;...\,;M\}\,$. Similarly, we can replace $\bk{\tilde{v}_l^t}{v_j^t}\, \bk{\tilde{v}_i^t}{v_k^t}$ with $\bk{\tilde{u}_l^t}{u_j^t}\, \bk{\tilde{v}_i^t}{v_k^t}$ or $\bk{\tilde{v}_l^t}{v_j^t}\, \bk{\tilde{u}_i^t}{u_k^t}\,$ for example.

\subsubsection{Deriving the System}
The system of partial differential equations is obtained by applying Itô's lemma to each of the three functions. Therefore, we detail how we obtain the equation on $S_V$ (the method for $S_U$ is almost identical) and the equation on $S_W\,$. We work with fixed $t$ and fixed $z\,,\tilde{z} \in \mathbb{C}\setminus \mathbb{R}\,$, all three independent of $M\,,N\,,m\,,n\,$.
\paragraph{First Equation}
Itô's formula on $S_V^{(N)}$ gives
\begin{align*}
dS_V^{(N)} &= \frac{1}{N} \, \sum_{i=1}^n \, \sum_{j=1}^N \, \frac{dV_{ij}(t)}{(\tilde{z} - \mu_i^t)(z - \lambda_j^t)} + \frac{1}{N} \, \sum_{i=1}^n \, \sum_{j=1}^N \, \frac{V_{ij}(t)}{(\tilde{z} - \mu_i^t)^2(z - \lambda_j^t)} \, d\mu_i^t \\
&\quad + \frac{1}{N} \, \sum_{i=1}^n \, \sum_{j=1}^N \, \frac{V_{ij}(t)}{(\tilde{z} - \mu_i^t)(z - \lambda_j^t)^2} \, d\lambda_j^t + \frac{1}{N} \, \sum_{i=1}^n \, \sum_{j=1}^N \, \frac{dV_{ij}(t)}{(\tilde{z} - \mu_i^t)^2(z - \lambda_j^t)}\, d\mu_i^t \\
&\quad + \frac{1}{N} \, \sum_{i=1}^n \, \sum_{j=1}^N \, \frac{dV_{ij}(t)}{(\tilde{z} - \mu_i^t)(z - \lambda_j^t)^2}\, d\lambda_j^t + \frac{1}{N} \, \sum_{i=1}^n \, \sum_{j=1}^N \, \frac{V_{ij}(t)}{(\tilde{z} - \mu_i^t)^2(z - \lambda_j^t)^2} \, d\mu_i^t \, d\lambda_j^t \\
&\quad + \frac{1}{N} \, \sum_{i=1}^n \, \sum_{j=1}^N \, \frac{V_{ij}(t)}{(\tilde{z} - \mu_i^t)^3(z - \lambda_j^t)} \, \left(d\mu_i^t\right)^2 + \frac{1}{N} \, \sum_{i=1}^n \, \sum_{j=1}^N \, \frac{V_{ij}(t)}{(\tilde{z} - \mu_i^t)(z - \lambda_j^t)^3} \, \left(d\lambda_j^t\right)^2 \,.
\end{align*}
Based on the correlations derived in Section \ref{sec:Dynamics}, we have:
\begin{itemize}
\item[$\bullet$] $dV_{ij}(t)\, d\mu_i^t =  \frac{4\, \sqrt{\mu_i^t}}{N}\, \sum\limits_{\substack{k = 1 \\ k \neq j}}^N \, \frac{\sqrt{\lambda_j^t}\, \ow{i}{i}{k}{j} + \sqrt{\lambda_k^t}\, \ow{i}{i}{j}{k}}{\lambda_j^t - \lambda_k^t}\, \ov{i}{i}{j}{k}\, dt = \mathcal{O}\left(\frac{1}{N^2}\right)\,.$
\item[$\bullet$] $dV_{ij}(t) \, d\lambda_j^t =  \frac{4\, \sqrt{\lambda_j^t}}{N} \, \sum\limits_{\substack{l=1\\ l \neq i}}^n \, \frac{\sqrt{\mu_i^t}\, \ow{l}{i}{j}{j} + \sqrt{\mu_l^t}\, \ow{i}{l}{j}{j}}{\mu_i^t - \mu_l^t} \, \ov{i}{l}{j}{j}\, dt = \mathcal{O}\left(\frac{1}{N^2}\right) \,.$
\item[$\bullet$] $d\mu_i^t\, d\lambda_j^t = \frac{4}{N} \, \sqrt{\mu_i^t \lambda_j^t}\, d\tilde{b}_i(t)\, db_j(t) = \frac{4}{N}\, \sqrt{\mu_i^t \lambda_j^t}\, W_{ij}(t) \, dt = \mathcal{O}\left(\frac{1}{N^2}\right)\,.$
\item[$\bullet$] $\left(d\mu_i^t\right)^2 = \frac{4}{N} \,\mu_i^t \, dt = \mathcal{O}\left(\frac{1}{N}\right)\,.$
\item[$\bullet$] $\left(d\lambda_j^t\right)^2 = \frac{4}{N} \, \lambda_j^t \, dt = \mathcal{O}\left(\frac{1}{N}\right)\,.$
\end{itemize}
Therefore, using the fact that $V_{ij}$ vanishes as $1\,/\,N$ in the scaling limit, we can rewrite our previous Itô formula as
\begin{align*}
dS_V^{(N)} &= \frac{1}{N} \, \sum_{i=1}^n \, \sum_{j=1}^N \, \frac{dV_{ij}(t)}{(\tilde{z} - \mu_i^t)(z - \lambda_j^t)} + \frac{1}{N} \, \sum_{i=1}^n \, \sum_{j=1}^N \, \frac{V_{ij}(t)}{(\tilde{z} - \mu_i^t)^2(z - \lambda_j^t)} \, d\mu_i^t \\
&\quad + \frac{1}{N} \, \sum_{i=1}^n \, \sum_{j=1}^N \, \frac{V_{ij}(t)}{(\tilde{z} - \mu_i^t)(z - \lambda_j^t)^2} \, d\lambda_j^t + o(1) \,.
\end{align*}
We denote the sums on the right-hand side respectively by $d\Sigma_V\,, d\Sigma_{\mu}$ and $d\Sigma_{\lambda}\,$. We can expand the first sum $d\Sigma_V$ using the dynamics of $V_{ij}\,$ for Appendix \ref{subsec:App_ito_overlaps} as
$$
d\Sigma_V = \left(I_{\mu} + I_{\lambda} + I_{\mu \lambda}\right)\, dt + dI_{\tilde{w}} + dI_w \,,
$$
where:
\begin{align*}
I_{\mu} &:= \frac{1}{N^2}\, \sum_{\substack{i,l = 1 \\ l \neq i}}^n \, \sum_{j=1}^N \, \frac{\left(\mu_i^t + \mu_l^t\right)\,\left(V_{lj} - V_{ij}\right)}{(\mu_i^t - \mu_l^t)^2(\tilde{z} - \mu_i^t)(z - \lambda_j^t)}\,,\\
I_{\lambda} &:= \frac{1}{N^2} \, \sum_{i=1}^n \, \sum_{\substack{j,k = 1 \\ k \neq j}}^N \, \frac{\left(\lambda_j^t + \lambda_k^t\right)\,\left(V_{ik} - V_{ij}\right)}{(\lambda_j^t - \lambda_k^t)^2(\tilde{z} - \mu_i^t)(z - \lambda_j^t)} \,, \\
I_{\mu \lambda} &:= \frac{2}{N^2} \, \sum_{\substack{i,l = 1 \\ l \neq i}}^n \, \sum_{\substack{j,k=1 \\ k \neq j}}^N \, \frac{\left[ V_{ij}\, \bar{W}_{lk} + \sqrt{\mu_i^t \lambda_j^t}\, \ov{i}{l}{k}{j}\, \ow{l}{i}{k}{j} \right]}{(\mu_i^t - \mu_l^t)(\lambda_j^t - \lambda_k^t)(\tilde{z} - \mu_i^t)(z - \lambda_j^t)} \,,\\
dI_{\tilde{w}} &:=  \frac{2}{N\,\sqrt{N}}\, \sum_{\substack{i,l = 1 \\ l \neq i}}^n \, \sum_{j=1}^N\, \frac{\sqrt{\mu_i^t}\, d\tilde{w}_{il}(t) + \sqrt{\mu_l^t}\, d\tilde{w}_{li}(t)}{(\mu_i^t - \mu_l^t)(\tilde{z} - \mu_i^t)(z - \lambda_j^t)}\, \ov{i}{l}{j}{j} \,, \\
dI_w &:= \frac{2}{N\,\sqrt{N}} \, \sum_{i=1}^n \, \sum_{\substack{j,k = 1 \\ k \neq j}}^N \, \frac{\sqrt{\lambda_j^t}\, dw_{jk}(t) + \sqrt{\lambda_k^t}\,dw_{kj}(t)}{(\lambda_j^t - \lambda_k^t)(\tilde{z} - \mu_i^t)(z - \lambda_j^t)} \, \ov{i}{i}{j}{k} \,.
\end{align*}
Our goal is to prove the following convergences:
\begin{itemize}
\item[$\bullet$] $I_{\mu}\, dt + d\Sigma_{\mu} \to \left(\alpha - \frac{\beta}{q} - 2 \alpha \tilde{z} \, \tilde{G}(\tilde{z},t)\right)\, \partial_{\tilde{z}}S_V \, dt - \alpha\,\tilde{G}(\tilde{z},t)\, S_V \, dt\,,$
\item[$\bullet$] $I_{\lambda} \, dt + d\Sigma_{\lambda} \to \left(1 - \frac{1}{q} - 2z \, G(z,t)\right) \, \partial_z S_V \, dt - G(z,t)\, S_V \, dt\,,$
\item[$\bullet$] $I_{\mu \lambda} \, dt \to 2\, S_W\, S_V\, dt\,,$
\item[$\bullet$] $dI_{\tilde{w}} \to 0$ and $dI_w \to 0\,.$
\end{itemize}
We start by manipulating $I_{\mu}\,$, applying symmetrisation (\ref{eq:S3}) to indices $i$ and $l$ with $a_{il} = \left(\mu_i^t + \mu_l^t\right)\, V_{lj} \,/\, (\mu_i^t - \mu_l^t)\,$, we transform it into
$$
I_{\mu} = \frac{1}{N^2} \, \sum_{\substack{i,l = 1 \\ l \neq i}}^n \, \sum_{j=1}^N \, \frac{\left(\mu_i^t + \mu_l^t\right)\, V_{lj}}{(\mu_i^t - \mu_l^t)(\tilde{z} - \mu_i^t)(\tilde{z} - \mu_l^t)(z - \lambda_j^t)} \,.
$$
Now, using the dynamics of $\mu_i^t\,$, we have
\begin{align*}
d\Sigma_{\mu} &= \frac{1}{N}\, \sum_{i=1}^n \, \sum_{j=1}^N \frac{V_{ij}}{(\tilde{z} - \mu_i^t)(z - \lambda_j^t)} \, \left( \frac{m}{N} \, dt + \frac{1}{N} \, \sum_{\substack{l = 1 \\ l \neq i}}^n \, \frac{\mu_i^t + \mu_l^t}{\mu_i^t - \mu_l^t} \, dt \right) \\
&\quad + \frac{2}{N\sqrt{N}} \, \sum_{i=1}^n\, \sum_{j=1}^N \, \frac{V_{ij}\, \sqrt{\mu_i^t}\, d\tilde{b}_i(t)}{(\tilde{z} - \mu_i^t)^2(z-\lambda_j^t)}\,,
\end{align*}
which, by inverting the indices $i$ and $l$ in the double sum, can be rewritten as
$$
d\Sigma_{\mu} = - \frac{m}{N} \, \partial_{\tilde{z}}S_V^{(N)}\, dt - \frac{1}{N^2} \, \sum_{\substack{i,l = 1 \\ l \neq i}}^n\, \sum_{j=1}^N\, \frac{\left(\mu_i^t + \mu_l^t\right)\, V_{lj}}{(\mu_i^t - \mu_l^t)(\tilde{z} - \mu_l^t)^2(z - \lambda_j^t)} \, dt + o(1) \,.
$$
Using identity (\ref{eq:I}), we obtain
$$
I_{\mu} \, dt + d\Sigma_{\mu} =\frac{1}{N^2} \, \sum_{\substack{i,l = 1 \\ l \neq i}}^n \, \sum_{j=1}^N
 \frac{\left(\mu_i^t + \mu_l^t\right)\, V_{lj}}{(\tilde{z} - \mu_l^t)^2(\tilde{z} - \mu_i^t)(z - \lambda_j^t)} \, dt - \frac{m}{N}\, \partial_{\tilde{z}}S_V^{(N)}\, dt + o(1) \,,$$
where we can add the diagonal terms $l = i$ (that are well defined since we got rid of the $\mu_i^t - \mu_l^t$ denominators) as their are vanishing in the scaling limit because of the factor $1\,/\,N^2$ and of the order of magnitude of $V_{lj}\,$. We can expand $(\mu_i^t + \mu_l^t)\,$ in the sum, the first sum we obtain is
\begin{align*}
\frac{1}{N^2}\, \sum_{\substack{i,l = 1}}^n \, \sum_{\substack{j = 1}}^N\,\frac{\mu_i^t \, V_{lj}}{(\tilde{z} - \mu_l^t)^2 (\tilde{z} - \mu_i^t)(z - \lambda_j^t)} &= \left(\frac{1}{N}\,\sum_{i = 1}^n\, \frac{\mu_i^t}{\tilde{z} - \mu_i^t}\right)\,\left( \frac{1}{N} \, \sum_{\substack{l =  1}}^n \, \sum_{j=1}^N\,\frac{V_{lj}}{(\tilde{z} - \mu_l^t)(z - \lambda_j^t)}\right) \\
&= \left(-\frac{n}{N} + \frac{n}{N}\tilde{z}\, \tilde{G}(\tilde{z},t)\right)\, \left(-\partial_{\tilde{z}}S_V^{(N)}\right) \\
&= \left(\frac{n}{N} - \frac{n}{N}\, \tilde{z}\, \tilde{G}(\tilde{z},t)\right)\, \partial_{\tilde{z}}S_V^{(N)} \,, 
\end{align*}
and the second one is
\begin{align*}
\frac{1}{N^2}\, \sum_{\substack{i,l = 1}}^n \, \sum_{j=1}^N\,\frac{\mu_l^t\, V_{lj}}{(\tilde{z} - \mu_l^t)^2 (\tilde{z} - \mu_i^t)(z - \lambda_j^t)} &= \left( \frac{1}{N} \, \sum_{i=1}^n\, \frac{1}{\tilde{z} - \mu_i^t} \right) \, \left( \frac{1}{N} \, \sum_{l=1}^n\, \sum_{j=1}^N\, \frac{\mu_l^t\, V_{lj}}{(\tilde{z} - \mu_l^t)^2(z-\lambda_j^t)} \right) \\
&= \frac{n}{N}\, \tilde{G}(\tilde{z},t) \left( -S_V^{(N)} - \tilde{z}\, \partial_{\tilde{z}}S_V^{(N)} \right) \,.
\end{align*}
Since $n\,/\,N \to \alpha$ and $m\,/\,N \to \beta \,/\, q\,$, we obtain the announced convergence
$$
I_{\mu}\, dt + d\Sigma_{\mu} \to \left(\alpha - \frac{\beta}{q} - 2 \alpha \tilde{z} \, \tilde{G}(\tilde{z},t)\right)\, \partial_{\tilde{z}}S_V - \alpha\,  \tilde{G}(\tilde{z},t)\, S_V \,.
$$
The method for the convergence of $I_{\lambda}\, dt + d\Sigma_{\lambda}$ is identical.

We now derive the limit of $I_{\mu \lambda}\,$, applying symmetrisation (\ref{eq:S4}) to it we obtain
$$
I_{\mu \lambda} = \frac{2}{N^2} \, \sum_{\substack{i,l = 1 \\ l \neq i}}^n\, \sum_{\substack{j,k = 1\\ k \neq j}}^N \, \frac{V_{ij}\, \bar{W}_{lk} + \sqrt{\mu_i^t \lambda_j^t}\, \ov{i}{l}{k}{j} \, \ow{l}{i}{k}{j}}{(\tilde{z} - \mu_i^t)(\tilde{z} - \mu_l^t)(z - \lambda_j^t)(z-\lambda_k^t)} \,.
$$
Expanding the numerator we get two sums $I_{\mu \lambda}^{(1)} + I_{\mu \lambda}^{(2)}\,$. Adding the diagonal terms $l = i$ and $k = j$ to the first one, since they vanish in the scaling limit, gives
\begin{align*}
I_{\mu \lambda}^{(1)} &= 2\, \left( \frac{1}{N} \, \sum_{\substack{i = 1}}^n \, \sum_{j=1}^N\, \frac{V_{ij}}{(\tilde{z} - \mu_i^t)(z - \lambda_j^t)}\right)\,\left(\frac{1}{N} \, \sum_{l=1}^n \, \sum_{k=1}^N \, \frac{ \bar{W}_{lk}}{(\tilde{z} - \mu_l^t)(z-\lambda_k^t)}\right) + o(1) \\
&= 2 \, S_W\, S_V + o(1)\,.
\end{align*}
Adding to the second sum its diagonal terms that are of order $1 \,/\, N\,$, we have
\begin{align*}
I_{\mu \lambda}^{(2)} &= \frac{2}{N^2}\, \sum_{i,l=1}^n \, \sum_{j,k=1}^N\, \frac{\sqrt{\mu_i^t \lambda_j^t} \, \bk{\tilde{v}_i^t}{v_k^t}\, \bk{\tilde{v}_l^t}{v_j^t}\,\bk{\tilde{v}_l^t}{v_k^t}\, \bk{\tilde{u}_i^t}{u_j^t}}{(\tilde{z} - \mu_i^t)(\tilde{z} - \mu_l^t)(z-\lambda_j^t)(z-\lambda_k^t)} + \mathcal{O}\left(\frac{1}{N}\right)
\end{align*}
Applying (\ref{eq:R}) with $c_{ij} = \sqrt{N}\,\sqrt{\mu_i^t \lambda_j^t} \, \bk{\tilde{u}_i^t}{u_j^t} = \mathcal{O}(1)\,$, we get $I_{\mu \lambda}^{(2)} \to 0\,$.

Finally, we prove that the Brownian terms $dI_{w}$ and $dI_{\tilde{w}}$ go to zero in the scaling limit. We detail the method for $dI_{w}$ only. The independence of $dw$ with respect to the other random variables in the sum indicates that $dI_w$ is centered. Furthermore, we can apply symmetrisation (\ref{eq:S3}) to the indices $j$ and $k$ which leads to 
$$
dI_w = \frac{2}{N\, \sqrt{N}}\, \sum_{i = 1}^{n}\, \sum_{\substack{j,k = 1}}^{N}\, \frac{\sqrt{\lambda_j^t}\, dw_{jk}(t)}{(\tilde{z} - \mu_i^t)(z - \lambda_j^t)(z - \lambda_k^t)}\, \ov{i}{i}{j}{k} \,,
$$
and we can write its variance as 
\begin{align*}
    \mathbb{E}\left[\left|dI_w  \right|^2 \right] &= \frac{4}{N^3} \, \mathbb{E} \left[\sum_{i,l = 1}^{n} \, \sum_{\substack{j,k,j',k' = 1 \\ k \neq j \\ k' \neq j'}}^{N} \frac{\sqrt{\lambda_j^t \lambda_{j'}^t} \, dw_{jk}(t) \, dw_{j'k'}(t)\, \ov{i}{i}{j}{k}\, \ov{l}{l}{j'}{k'}}{(\tilde{z} - \mu_i^t)(\tilde{z}^* - \mu_l^t)(z - \lambda_j^t)(z - \lambda_k^t)(z^*- \lambda_{j'}^t)(z^* - \lambda_{k'}^t)}  \right] \\
    &= \frac{4}{N^3}\, \mathbb{E}\left[\sum_{i,l = 1}^{n} \, \sum_{\substack{j,k = 1 \\ k \neq j}}^{N} \, \frac{\lambda_j^t \, \bk{\tilde{v}_i^t}{v_j^t}\, \bk{\tilde{v}_i^t}{v_k^t}\, \bk{\tilde{v}_l^t}{v_j^t}\, \bk{\tilde{v}_l^t}{v_k^t}}{(\tilde{z} - \mu_i^t)(\tilde{z}^* - \mu_l^t)\left| z - \lambda_j^t\right|^2 \left|z - \lambda_k^t\right|^2}\right]\, dt \,,
\end{align*}
which gives $\mathbb{E}\left[\left|dI_w\right|^2\right] = \mathcal{O}(1 \,/\, N^2)$ using (\ref{eq:R}) with $c_{ij} = \sqrt{N}\,\lambda_j^t \, \bk{\tilde{v}_i^t}{v_j^t} = \mathcal{O}(1)\,$. Since the variances are summable with respect to $N\,$, Borel-Cantelli's lemma indicates that $dI_w \to 0$ almost surely.

We have proved that randomness vanishes almost surely in the equation on $S_V^{(N)}$ and leads to 
$$
\partial_t S_V = g(z,t)\, \partial_z S_V + \tilde{g}(\tilde{z}, t)\, \partial_{\tilde{z}} S_V + \left(2 \, S_W - G(z,t) - \alpha \, \tilde{G}(\tilde{z},t)\right) \, S_V \,,
$$
with $g(z,t) := 1 - \frac{1}{q} - 2z \, G(z,t)$ and $\tilde{g}(\tilde{z}, t) := \alpha - \frac{\beta}{q} - 2 \alpha \tilde{z}\, \tilde{G}(\tilde{z},t)\,$.

\paragraph{Second Equation}
The equation on $S_U$ is obtained with the same method. The only difference comes from the fact that instead of obtaining the term $G_N(z,t) \, S_U^{(N)}\,$, we get
$$
\frac{1}{N} \, \sum_{j = 1}^{M} \, \frac{1}{z - \lambda_j^t} \, S_U^{(N)} \,,
$$
which is equal to (recalling that we introduced the notations $\lambda_{N+1}^t = ... = \lambda_M^t = 0$ for simplicity)
$$
\left(\frac{M - N}{N\, z} + G_N(z,t) \right) \, S_U^{(N)}\,,
$$
that converges to
$$
\left(\frac{\frac{1}{q} - 1}{z} + G(z,t)\right)\, S_U \,.
$$
A similar modification is obtained for the $\tilde{G}(\tilde{z},t) \, S_U$ term.

\paragraph{Third Equation}
For the equation on $S_W\,$, once summed, the Brownian terms almost surely vanish in the scaling limit using the same argument as for $S_V\,$. Likewise, the sums of the form 
$$
\frac{1}{N^2}\, \sum_{\substack{i,l = 1 \\ l \neq i}}\, \sum_{\substack{j,k = 1\\ k \neq j}}\, \frac{c_{ij}\, \left[\ov{l}{l}{j}{k}\, \ou{i}{i}{j}{k}\right]}{(\tilde{z} - \mu_i^t)(\mu_i^t - \mu_l^t)(z - \lambda_j^t)(\lambda_j^t - \lambda_k^t)}
$$
go to zero (using arguments similar to (\ref{eq:R})). Therefore, we focus on the transformation of the non vanishing terms (we recall that in Appendix \ref{subsec:App_ito_overlaps} we introduced the notation $\bar{W}_{ij} = \sqrt{\mu_i^t \lambda_j^t}\, W_{ij}$):
\begin{align*}
    dS_W^{(N)} &= \frac{1}{N}\, \sum_{i=1}^n\, \sum_{j=1}^{N}\, \frac{\sqrt{\mu_i^t \lambda_j^t}\, dW_{ij}(t)}{(\tilde{z} - \mu_i^t)(z - \lambda_j^t)} \\
    &\quad + \frac{1}{N} \, \sum_{i=1}^{n}\, \sum_{j = 1}^{N}\,\left(\frac{\bar{W}_{ij}(t)}{2\mu_i^t(\tilde{z} - \mu_i^t)(z - \lambda_j^t)} + \frac{\bar{W}_{ij}(t)}{(\tilde{z} - \mu_i^t)^2(z - \lambda_j^t)}\right)\, d\mu_i^t \\
    &\quad + \frac{1}{N}\, \sum_{i= 1}^{n}\, \sum_{j=1}^{N}\, \left(\frac{\bar{W}_{ij}(t)}{2 \lambda_j^t(\tilde{z} - \mu_i^t)(z- \lambda_j^t)} + \frac{\bar{W}_{ij}(t)}{(\tilde{z} - \mu_i^t)(z - \lambda_j^t)^2}\right)\, d\lambda_j^t + o(1)\,.
\end{align*}
We denote by $\Sigma_W\,$, $\Sigma_{\mu}$ and $\Sigma_{\lambda}$ the three sums on the right hand side, in their respective order. From what we said, most of the terms vanish in $\Sigma_W$ so we can write
$$
\Sigma_W = \left(I_{\mu} + I_{\lambda} + I_{VU} + I_{W}\right)\, dt + o(1)\,,
$$
where:
\begin{itemize}
    \item[$\bullet$] $I_{\mu} := \frac{1}{N^2}\, \sum\limits_{\substack{i,l = 1 \\ l \neq i}}^{n}\, \sum\limits_{j=1}^{N}\, \frac{2\mu_i^t\, \bar{W}_{lj} - (\mu_i^t + \mu_l^t)\, \bar{W}_{ij}}{(\mu_i^t - \mu_l^t)^2(\tilde{z} - \mu_i^t)(z - \lambda_j^t)} + \frac{n-m}{2N^2}\, \sum\limits_{i=1}^{n}\, \sum\limits_{j=1}^{N}\, \frac{\bar{W}_{ij}}{\mu_i^t(\tilde{z} - \mu_i^t)(z - \lambda_j^t)} \,, $ 
    \item[$\bullet$] $I_{\lambda} := \frac{1}{N^2}\, \sum\limits_{i=1}^{n}\, \sum\limits_{\substack{j,k=1 \\ k \neq j}}^{N}\, \frac{2\lambda_j^t\, \bar{W}_{ik} - (\lambda_j^t + \lambda_k^t)\, \bar{W}_{ij}}{(\lambda_j^t - \lambda_k^t)^2(\tilde{z} - \mu_i^t)(z - \lambda_j^t)} + \frac{N-M}{2N^2}\, \sum\limits_{i=1}^{n}\, \sum\limits_{j=1}^{N}\, \frac{\bar{W}_{ij}}{\lambda_j^t(\tilde{z} - \mu_i^t)(z - \lambda_j^t)} \,,$
    \item[$\bullet$] $I_{W} := \frac{1}{N^2} \, \sum\limits_{\substack{i,l = 1 \\ l \neq i}}^{n}\, \sum\limits_{\substack{j,k = 1 \\ k \neq j}}^N \, \frac{\sqrt{\mu_i^t \mu_l^t \lambda_j^t \lambda_k^t}\, \left[W_{ij}\,W_{lk}\right]}{(\mu_i^t - \mu_l^t)(\tilde{z} - \mu_i^t)(\lambda_j^t - \lambda_k^t)(z - \lambda_j^t)} \,,$
\end{itemize}
and
\begin{align*}
    I_{VU} :&= \frac{1}{N^2} \, \sum\limits_{\substack{i,l = 1 \\ l \neq i}}^{n}\, \sum\limits_{\substack{j,k = 1 \\ k \neq j}}^N \, \frac{\mu_i^t \lambda_j^t\, \left[V_{ij}\,U_{lk}\right]}{(\mu_i^t - \mu_l^t)(\tilde{z} - \mu_i^t)(\lambda_j^t - \lambda_k^t)(z - \lambda_j^t)} \\
    &\quad + \frac{1}{N^2}\, \sum_{i=1}^n \, \sum_{l = n+1}^m \, \sum_{\substack{j,k=1 \\ k \neq j}}^{N}\, \frac{\lambda_j^t \, \left(V_{ij}\, U_{lk} + V_{ik}\, U_{lj}\right)}{(\tilde{z} - \mu_i^t)(\lambda_j^t - \lambda_k^t)(z - \lambda_j^t)}
\end{align*}
\begin{align*}
    &\quad + \frac{1}{N^2}\, \sum_{\substack{i,l=1 \\ l \neq i}}^{n}\, \sum_{j=1}^{N}\, \sum_{k = N+1}^{M}\, \frac{\mu_i^t\, \left(V_{ij}\, U_{lk} + V_{lj}\,U_{ik}\right)}{(\mu_i^t - \mu_l^t)(\tilde{z} - \mu_i^t)(z - \lambda_j^t)} \\
    &\quad + \frac{1}{N^2}\, \sum_{i=1}^{n}\, \sum_{l = n+1}^{m}\, \sum_{j=1}^{N}\, \sum_{k = N+1}^{M}\, \frac{V_{ij}\, U_{lk}}{(\tilde{z} - \mu_i^t)(z - \lambda_j^t)}\,.
\end{align*}
We are going to prove the following convergences:
\begin{itemize}
    \item[$\bullet$] $I_{\mu}\, dt +\Sigma_{\mu} \to \left(\alpha - \frac{\beta}{q} - 2 \alpha \tilde{z}\, \tilde{G}(\tilde{z},t)\right)\, \partial_{\tilde{z}}S_W\,dt\,,$
    \item[$\bullet$] $I_{\lambda}\, dt + \Sigma_{\lambda} \to \left(1 - \frac{1}{q} - 2z \, G(z,t)\right)\, \partial_z S_W \, dt \,,$
    \item[$\bullet$] $I_{VU} \to z \tilde{z}\, S_V \, S_U\,,$
    \item[$\bullet$] $I_W \to S_W^2\,.$   
\end{itemize}
We begin with the convergence of $I_{\mu}\, dt + \Sigma_{\mu}\,$. We can rewrite the first sum in $I_{\mu}$ as 
$$
\frac{1}{N^2}\, \sum_{\substack{i,l = 1 \\ l \neq i}}^{n}\, \sum_{j=1}^{N}\, \frac{\mu_i^t\,\left(\bar{W}_{lj} - \bar{W}_{ij}\right)}{(\mu_i^t - \mu_l^t)^2(\tilde{z} - \mu_i^t)(z - \lambda_j^t)} + \frac{1}{N^2}\, \sum_{\substack{i,l = 1 \\l \neq i}}^{n}\, \sum_{j=1}^{N}\, \frac{\mu_i^t\, \bar{W}_{lj} - \mu_l^t \, \bar{W}_{ij}}{(\mu_i^t - \mu_l^t)^2(\tilde{z}-\mu_i^t)(z-\lambda_j^t)}\,.
$$
In the first term, we can replace $\mu_i^t$ in the numerator by $\tilde{z}$ because the difference between the two sums is a null sum (the summand is antisymmetric with respect to $i$ and $l$). Applying symmetrisation (\ref{eq:S3}) to both sums, we obtain
$$
\frac{1}{N^2}\, \sum_{\substack{i,l=1\\l \neq i}}^{n}\, \sum_{j=1}^N\, \frac{(\tilde{z} + \mu_i^t)\, \bar{W}_{lj}}{(\mu_i^t - \mu_l^t)(\tilde{z} - \mu_i^t)(z - \lambda_j^t)}\,.
$$
Once again, we use identity (\ref{eq:I}) to transform it into
$$
\frac{1}{N^2} \, \sum_{\substack{i,l=1\\l \neq i}}^{n}\, \sum_{j=1}^{N}\, \frac{(\tilde{z} + \mu_i^t)\, \bar{W}_{lj}}{(\tilde{z} - \mu_l^t)^2(\mu_i^t - \mu_l^t)(z - \lambda_j^t)} + \frac{1}{N^2} \, \sum_{\substack{i,l=1\\l \neq i}}^{n}\, \sum_{j=1}^{N}\, \frac{(\tilde{z} + \mu_i^t)\, \bar{W}_{lj}}{(\tilde{z} - \mu_l^t)^2(\tilde{z} - \mu_i^t)(z - \lambda_j^t)}\,.
$$
The second sum converges to 
$$
\left(\alpha - 2 \alpha \tilde{z}\,\tilde{G}(\tilde{z}, t)\right)\, \partial_{\tilde{z}}S_W\,,
$$
and we denote by $A$ the first sum that we will combine with $\Sigma_{\mu}\,$. We recall that 
$$
d\mu_i^t = \frac{m}{N}\, dt + \frac{1}{N}\, \sum_{\substack{l=1 \\ l \neq i}}^{n}\, \frac{\mu_i^t + \mu_l^t}{\mu_i^t - \mu_l^t}\, dt + o(1)\,,
$$
so that
\begin{align*}
    \Sigma_{\mu} &= \frac{m}{2N^2}\, \sum_{i=1}^{n}\, \sum_{j=1}^N\, \frac{\bar{W}_{ij}}{\mu_i^t (\tilde{z} - \mu_i^t)(z - \lambda_j^t)}\, dt + \frac{m}{N^2}\, \sum_{i=1}^n\, \sum_{j=1}^{N}\, \frac{\bar{W}_{ij}}{(\tilde{z} - \mu_i^t)^2(z - \lambda_j^t)}\, dt \\
    &\quad + \frac{1}{N^2}\, \sum_{\substack{i,l=1 \\ l \neq i}}^{n}\, \sum_{j=1}^N\, \frac{\mu_i^t+ \mu_l^t}{\mu_i^t - \mu_l^t}\, \left(\frac{\bar{W}_{ij}}{2\mu_i^t (\tilde{z} - \mu_i^t)(z - \lambda_j^t)} + \frac{\bar{W}_{ij}}{(\tilde{z} - \mu_i^t)^2(z - \lambda_j^t)}\right)\, dt + o(1)\,,
\end{align*}
where the second sum converges to $-\frac{\beta}{q}\, \partial_{\tilde{z}}S_W\, dt$ and the last sum, if added to $A\,dt$ (after exchanging the indices $i$ and $l$), equals 
$$
-\frac{n}{N^2}\, \sum_{i=1}^{n}\, \sum_{j=1}^{N}\, \frac{\bar{W}_{ij}}{2\mu_i^t(\tilde{z} - \mu_i^t)(z - \lambda_j^t)}\,dt\,.
$$
Therefore,
\begin{align*}
    \Sigma_{\mu} + A\,dt &= \frac{m-n}{2N^2}\,\sum_{i=1}^{n}\, \sum_{j=1}^{N}\, \frac{\bar{W}_{ij}}{\mu_i^t(\tilde{z} - \mu_i^t)(z - \lambda_j^t)}\, dt -\frac{\beta}{q}\, \partial_{\tilde{z}}S_W\, dt + o(1)
\end{align*}
which cancels out with the second sum in the definition of $I_{\mu}\,dt\,$. Finally, we have proved that 
$$
I_{\mu} \, dt + \Sigma_{\mu} \longrightarrow \left(\alpha - \frac{\beta}{q} - 2\alpha \tilde{z}\, \tilde{G}(\tilde{z}, t)\right)\, \partial_{\tilde{z}}S_W\, dt\,.
$$
The demonstration for the convergence of $I_{\lambda}\, dt + \Sigma_{\lambda}$ is identical.

For $I_W\,$, we first notice that $\sqrt{\mu_i^t \mu_l^t \lambda_j^t \lambda_k^t}\, \left[W_{ij}\, W_{lk}\right] = \left[\bar{W}_{ij}\, \bar{W}_{lk}\right]\,$. Then, applying symmetrisation (\ref{eq:S4}) we get 
$$
I_W = \frac{1}{N^2}\, \sum_{\substack{i,l=1\\ l \neq i}}^{n}\, \sum_{\substack{j,k = 1 \\ k \neq j}}^{N}\, \frac{\bar{W}_{ij}\, \bar{W}_{lk}}{(\tilde{z} - \mu_i^t)(\tilde{z} - \mu_l^t)(z - \lambda_j^t)(z - \lambda_k^t)} \,,
$$
which converges to $S_W^2\,$.

We now focus on the remaining term $I_{VU}\,$. Considering its first sum, one can write $\mu_i^t \lambda_j^t = (\mu_i^t - \tilde{z})\, \lambda_j^t + (\lambda_j^t - z)\, \tilde{z} + z \, \tilde{z}$ to see that we can replace $\mu_i^t \lambda_j^t$ by $z\tilde{z}$ in the numerator because the difference between the two sums are two null sums (antisymmetric with respect to $i$ and $l$ or to $j$ and $k$). Therefore, the first sum in $I_{VU}$ equals, after applying symmetrisation (\ref{eq:S4}),
$$
\frac{z\tilde{z}}{N^2}\, \sum_{\substack{i,l = 1 \\ l \neq i}}^n \sum_{\substack{j,k=1 \\ k \neq j}}^{N}\, \frac{V_{ij}\, U_{lk}}{(\tilde{z} - \mu_i^t)(\tilde{z} - \mu_l^t)(z - \lambda_j^t)(z - \lambda_k^t)}\,.
$$
The same type of reasoning can be applied to the other sums composing $I_{VU}$ until we obtain
\begin{align*}
    I_{VU} &= z\tilde{z}\, S_V^{(N)} \, \left(\frac{1}{N}\, \sum_{l = 1}^n \, \sum_{k = 1}^{N}\, \frac{U_{lk}}{(\tilde{z} - \mu_l^t)(z - \lambda_k^t)} + \frac{1}{N}\, \sum_{l=n+1}^{m}\, \sum_{k = 1}^{N}\, \frac{U_{lk}}{\tilde{z} (z - \lambda_k^t)}\right. \\
    &\quad \left.+ \frac{1}{N}\, \sum_{l=1}^{n}\, \sum_{k=N+1}^{M}\, \frac{U_{lk}}{(\tilde{z} - \mu_l^t)z}+ \frac{1}{N}\, \sum_{l=n+1}^{m}\, \sum_{k=N+1}^{M}\, \frac{U_{lk}}{\tilde{z} z}\right) + o(1) \\
    &= z\tilde{z}\, S_V^{(N)} \, S_U^{(N)} + o(1)\,.
\end{align*}
Thus, $I_{VU}$ converges to $z\tilde{z}\, S_V \, S_U\,$.

Finally, we have proven that $S_W$ satisfies the announced deterministic differential equation.

\subsection{Solving the System}
\label{subsec:solving_system}
Let $z\,, \tilde{z} \in \mathbb{C} \setminus \mathbb{R}$ and $t\geq 0\,$. We introduce a new variable $s\,$, as well as functions $z(s)\,$, $\tilde{z}(s)$ and $t(s)\,$ such that $z(0) = z\,$, $\tilde{z}(0) = \tilde{z}$ and $t(0) = t\,$. Moreover, we introduce the notation $\hat{S_V}(s) := S_V(z(s)\,,\tilde{z}(s)\,, t(s))$ and similarly for our other functions in the equations. Denoting by $c$ (respectively $\tilde{c}$) the constant $\frac{1}{q} - 1$ (respectively $\frac{\beta}{q} - \alpha$), if 
$$
\begin{cases}
    t'(s) = 1\\
    z'(s) = 2 \hat{G}(s)\, z(s) + c \\
    \tilde{z}'(s) = 2 \alpha \hat{\tilde{G}}(s)\, \tilde{z}(s) + \tilde{c}\,,
\end{cases}
$$
then the chain rule gives
$$
\begin{cases}
    \hat{S_V}'(s) = \left(2 \, \hat{S_W}(s) - \hat{G}(s) - \alpha \, \hat{\tilde{G}}(s)\right) \, \hat{S_V}(s) \\
    \hat{S_U}'(s) = \left(2 \, \hat{S_U}(s) - \frac{c}{z(s)} - \hat{G}(s) - \frac{\tilde{c}}{\tilde{z}(s)} - \alpha \, \hat{\tilde{G}}(s)\right) \, \hat{S_V}(s) \\
    \hat{S_W}'(s) = \hat{S_W}^2(s) + z(s) \tilde{z}(s)\, \hat{S_V}(s)\, \hat{S_U}(s)\,.
\end{cases}
$$
Additionally, under the previous conditions on $z(s)\,$, $\tilde{z}(s)$ and $t(s)$ we know from equation (\ref{eq:G_Burgers}) and its resolution in Appendix \ref{subsec:App_burgers} that \begin{itemize}
    \item[$\bullet$] $t(s) = t + s\,,$
    \item[$\bullet$] $z(s) = \left(1 + s\,\hat{G}(0)\right)\, \left(z \,(1 + s \hat{G}(0)) + cs\right)\,,$
    \item[$\bullet$] $\hat{G}(s) = \frac{\hat{G}(0)}{1 + s \, \hat{G}(0)}\,.$
\end{itemize}
The equation (\ref{eq:G_tilde_Burgers}) on $\tilde{G}$ can give us similarly:
\begin{itemize}
    \item[$\bullet$] $\tilde{z}(s) = \left(1 + \alpha s \, \hat{\tilde{G}}(0)\right)\, \left(\tilde{z}\, (1 + \alpha s \, \hat{\tilde{G}}(0)) + \tilde{c}s\right)\,,$
    \item[$\bullet$] $\hat{\tilde{G}}(s) = \frac{\hat{\tilde{G}}(0)}{1 + \alpha s \, \hat{\tilde{G}}(0)}\,.$
\end{itemize}
Therefore, the equations on $\hat{S_V}$ and $\hat{S_U}$ lead, after integration, to
$$
\hat{S_V}(s) = \frac{\hat{S_V}(0)}{\left(1 + s\,\hat{G}(0)\right)\,\left(1 + \alpha s \, \hat{\tilde{G}}(0)\right)}\, e^{2\, \int_0^s\, \hat{S_W}(u)\, du}\,,
$$
$$
\hat{S_U}(s) = \frac{z\tilde{z}\, \hat{S_U}(0)}{\left(z \,(1 + s \hat{G}(0)) + cs\right)\,\left(\tilde{z}\, (1 + \alpha s \, \hat{\tilde{G}}(0)) + \tilde{c}s\right)}\, e^{2\, \int_0^s\, \hat{S_W}(u)\, du}\,.
$$
This leaves us with the following differential equation on $\hat{S_W}\,$:
$$
\hat{S_W}'(s) = \hat{S_W}^2(s) + z\tilde{z}\, \hat{S_V}(0)\, \hat{S_U}(0)\, e^{4\,\int_0^s\, \hat{S_W}(u)\, du} \,.
$$
We are going to solve it explicitly. For readability we introduce the notations $f := \hat{S_W}\,$, $F := \int_0^{\cdot}\,f(u)\, du$ and $a := z\tilde{z}\, \hat{S_V}(0)\, \hat{S_U}(0)\,$. With the change of variable $x = F(s)\,$, we get
$$
\frac{df}{dx}\, f = f^2 + a\, e^{4x}\,,
$$
so that $g := f^2$ satisfies
$$
\frac{dg}{dx} = 2\,g + 2a\,e^{4x}\,,
$$
which gives, 
$$
g(x) = (f^2(0) - a)\, e^{2x} + a\, e^{4x}\,.
$$
This can be rewritten into an order 1 differential equation on $F\,$,
$$
\frac{dF}{ds} = \pm \sqrt{(f^2(0) - a)\, e^{2F} + a\,e^{4F}}\,.
$$
We separate the variables and integrate, which leads to 
$$
\sqrt{a + (f^2(0) - a)\,e^{-2F(s)}} = \sqrt{f^2(0)} \pm (f^2(0) - a)\, s \,,
$$
and finally,
$$
F(s) = - \frac{1}{2} \, \log\left(1 \pm 2 \, \sqrt{f^2(0)}\,s + (f^2(0) - a)\, s^2\right) \,.
$$
We can now differentiate to obtain
$$
f(s) = \frac{\mp \sqrt{f^2(0)} - (f^2(0) - a)\, s}{1 \pm 2 \, \sqrt{f^2(0)}\,s + (f^2(0) - a)\, s^2}\,.
$$
The condition at $s=0$ gives $f(0) = \mp \sqrt{f^2(0)}\,$, therefore we end up with
$$
f(s) = \frac{f(0) + (a - f^2(0))\, s}{1 - 2 \, f(0)\,s - (a - f^2(0))\, s^2}\,.
$$
Putting all of this together, we obtain the system
$$
\begin{cases}
    \hat{S_W}(s) = \frac{\hat{S_W}(0) + \left(z\tilde{z}\, \hat{S_V}(0) \, \hat{S_U}(0) - \hat{S_W}^2(0)\right)\, s}{1 - 2 \, \hat{S_W}(0)\, s - \left(z\tilde{z}\, \hat{S_V}(0) \, \hat{S_U}(0) - \hat{S_W}^2(0)\right)\, s^2} \\
    \hat{S_V}(s) = \frac{\hat{S_V}(0)}{\left(1 + s \hat{G}(0)\right)\, \left(1 + \alpha s \, \hat{\tilde{G}}(0)\right)\, \left(1 - 2 \, \hat{S_W}(0)\, s - \left(z\tilde{z}\, \hat{S_V}(0) \, \hat{S_U}(0) - \hat{S_W}^2(0)\right)\, s^2\right)}\\
    \hat{S_U}(s) = \frac{z\tilde{z}\,\hat{S_U}(0)}{\left(z\, (1 + s \hat{G}(0)) + cs\right)\, \left(\tilde{z}\, (1 + \alpha s \, \hat{\tilde{G}}(0)) + \tilde{c}s\right)\, \left(1 - 2 \, \hat{S_W}(0)\, s - \left(z\tilde{z}\, \hat{S_V}(0) \, \hat{S_U}(0) - \hat{S_W}^2(0)\right)\, s^2\right)} \,.
\end{cases}
$$
We denote by $D(s)$ the common denominator 
$$\left(1 - 2 \, \hat{S_W}(0)\, s - \left(z\tilde{z}\, \hat{S_V}(0) \, \hat{S_U}(0) - \hat{S_W}^2(0)\right)\, s^2\right)^{-1}\,.
$$
One can solve for $D$ using the previous system of equations, which gives
$$
D(s) = (1 + \hat{S_W}(s)\, s)^2 - z(s)\tilde{z}(s)\, \hat{S_V}(s) \, \hat{S_U}(s) \, s^2 \,.
$$
Thus, we can invert the system:
$$
\begin{cases}
    \hat{S_W}(0) = \frac{\hat{S_W}(s)\, (1 + \hat{S_W}(s)\,s) - z(s) \tilde{z}(s) \, \hat{S_V}(s)\, \hat{S_U}(s) s}{D(s)} \\
    \hat{S_V}(0) = \frac{\left(1 + s \,\hat{G}(0)\right)\, \left(1 + \alpha s \, \hat{\tilde{G}}(0)\right)\, \hat{S_V}(s)}{D(s)}\\
    \hat{S_U}(s) = \frac{\left(z\, (1 + s \,\hat{G}(0)) + cs\right)\, \left(\tilde{z}\, (1 + \alpha s \, \hat{\tilde{G}}(0)) + \tilde{c}s\right)\, \hat{S_U}(s)}{z\tilde{z}\, D(s)} \,.
\end{cases}
$$
Finally, since $\hat{f}(0) = f(z\,, \tilde{z}\,, t)\,$, evaluating at $s = -t$ gives the announced result.

\subsection{Inversion in the Marchenko-Pastur Case}
\label{subsec:inversion_MP}
We detail the case of $\bar{V}$ as the other functions are obtained almost identically.
First, we recall that 
$$
\lim_{\varepsilon \to 0^+}\, G(\lambda \pm i \, \varepsilon, t) = v(\lambda, t) \mp i \, \pi \, \rho(\lambda, t)
$$
where 
$$
\rho(\lambda, t) = \frac{\sqrt{\left((1 + \frac{1}{\sqrt{q}})^2\, t - \lambda\right)\, \left(\lambda - (1 - \frac{1}{\sqrt{q}})^2\, t\right)}}{2\pi \lambda t}
$$
and 
$$
v(\lambda , t) = \frac{\lambda - (\frac{1}{q} - 1)\, t}{2 \lambda t} \,.
$$
We have a similar relation between $\tilde{G}$ and 
$$
\tilde{\rho}(\mu, t) = = \frac{\sqrt{\left((\sqrt{\alpha} + \sqrt{\frac{\beta}{q}})^2\, t - \mu\right)\, \left(\mu - (\sqrt{\alpha} - \sqrt{\frac{\beta}{q}})^2\, t\right)}}{2\pi \alpha \mu t}\,,
$$
$$
\tilde{v}(\mu, t) = \frac{\mu - (\frac{\beta}{q} - \alpha)\,t}{2 \alpha \mu t} \,.
$$
Therefore if we define $S_V^{\pm} := \lim_{\varepsilon \to 0^+}\, S_V(\lambda - i \, \varepsilon, \mu \pm i \, \varepsilon, t)\,$, then,
$$
S_V^{\pm} = \frac{\alpha \, \left(A - i B\right)\, \left(\tilde{A} \pm i \tilde{B}\right)}{\left(A - i B\right)\,\left(\lambda\, A - ct - i \, \lambda \, B\right) \left(\tilde{A} \pm i \tilde{B}\right)\, \left(\mu \, \tilde{A} - \tilde{c}t \pm i \, \mu \, \tilde{B}\right) - \frac{\alpha \beta}{q}\, t^2} \,,
$$
where:\begin{itemize}
    \item[$\bullet$] $A := 1 - t\,v(\lambda,t)\,,$
    \item[$\bullet$] $B := \pi t \, \rho(\lambda,t)\,,$
    \item[$\bullet$] $\tilde{A} := 1 - \alpha t \, \tilde{v}(\mu,t)\,,$
    \item[$\bullet$] $\tilde{B} := \alpha \pi t \, \tilde{\rho}(\mu,t)\,,$
    \item[$\bullet$] $c := \frac{1}{q} - 1\,,$
    \item[$\bullet$] $\tilde{c} := \frac{\beta}{q} - \alpha\,.$     
\end{itemize}
We can simplify this into
$$
S_V^{\pm} = \frac{\alpha \, \left(A - i \, B\right)\, \left(\tilde{A} \pm i \, \tilde{B}\right)}{\left(A\, (\lambda\, A - ct) - \lambda \, B^2 - i\, B \, (2\lambda \, A - ct)\right)\, \left(\tilde{A}\, (\mu \, \tilde{A} - \tilde{c}t) - \mu \, \tilde{B}^2 \pm i \, \tilde{B}\, (2\mu \, \tilde{A} - \tilde{c}t)\right) - \frac{\alpha \beta}{q}\, t^2}\,.
$$
This form is very practical since we remark that $A = \frac{\lambda + ct}{2 \lambda}$ and $\tilde{A} = \frac{\mu + \tilde{c} t}{2 \mu}\,$, therefore $2\lambda \, A - ct = \lambda$ and $2 \mu \, \tilde{A} - \tilde{c}t = \mu\,$. Furthermore, rewriting $B$ and $\tilde{B}$ leads to 
$$
B^2 = \frac{-\lambda^2 + 2\,(1 + \frac{1}{q})\, t \, \lambda - c^2 t^2}{4 \,\lambda^2} \quad \text{and} \quad \tilde{B}^2 = \frac{-\mu^2 + 2 \, (\alpha + \frac{\beta}{q})\, t \, \mu - \tilde{c}^2 t^2}{4 \,\mu^2} \,.
$$
Therefore, $A\, (\lambda \, A - ct) - \lambda\,B^2 = \frac{\bar{\lambda}}{2}$ where $\bar{\lambda} := \lambda - \left( 1 + \frac{1}{q}\right)\, t$ and similarly $\tilde{A}\, (\mu \, \tilde{A} - \tilde{c}t) - \mu\, \tilde{B}^2 = \frac{\bar{\mu}}{2}$ where $\bar{\mu} := \mu - \left(\alpha + \frac{\beta}{q}\right)\, t\,$. We end up with
$$
S_V^{\pm} = \frac{\alpha \, \left(A - i \, B\right)\, \left(\tilde{A} \pm i \, \tilde{B}\right)}{\left(\frac{\bar{\lambda}}{2} - i \, \lambda \, B\right)\, \left(\frac{\bar{\mu}}{2} \pm i \, \mu \tilde{B}\right) - \frac{\alpha \beta}{q}\, t^2} =: \frac{N_{\pm}}{D_{\pm}} \,.
$$
In order to compute $\bar{V}\,$, we need to explicit the real part of $S_V^+ - S_V^-\,$. We have
\begin{align*}
    S_V^+ - S_V^- &= \frac{N_+\, D_- - N_- \, D_+}{D_+ \, D_-} \\
    &= \frac{\left(N_+\, D_- - N_- \, D_+\right)\, D_+^*\, D_-^*}{\left|D_+\, D_- \right|^2} \,,
\end{align*}
so we begin with simplifying the denominator. When needed, we use the fact that $\lambda^2 \, B^2 = \frac{t^2}{q} - \frac{\bar{\lambda}^2}{4}$ and $\mu^2 \, \tilde{B}^2 = \frac{\alpha \beta}{q}\, t^2 - \frac{\bar{\mu}^2}{4}\,$.
\begin{align*}
    \left| D_+ \, D_-\right|^2 &= \left| \left(\frac{\bar{\lambda}}{2} - i\, \lambda \, B\right)^2\, \left(\frac{\bar{\mu}^2}{4} + \mu^2\, \tilde{B}^2\right) - \frac{\alpha \beta}{q} \, t^2 \, \left(\frac{\bar{\lambda}}{2} - i\, \lambda \, B\right)\, \bar{\mu} + \frac{\alpha^2 \beta^2}{q^2}\, t^4 \right|^2 \\
    &= \left|\left(\frac{\bar{\lambda}^2}{4} - \lambda^2\, B^2 - i \, \lambda \, \bar{\lambda}\, B\right)\, \frac{\alpha \beta}{q}\, t^2 - \frac{\alpha \beta}{q} \, t^2 \, \left(\frac{\bar{\lambda}}{2} - i\, \lambda \, B\right)\, \bar{\mu} + \frac{\alpha^2 \beta^2}{q^2}\, t^4  \right|^2 \\
    &= \frac{\alpha^2 \beta^2}{q^2}\, t^4 \,\left| \frac{\bar{\lambda}^2}{2} - \frac{t^2}{q} - \frac{\bar{\lambda}\, \bar{\mu}}{2} + \frac{\alpha \beta}{q}\, t^2 - i \, \lambda \, \bar{\lambda}\, B + i \, \lambda\, \bar{\mu} \, B\right|^2 \\
    &= \frac{\alpha^2 \beta^2}{q^2}\, t^4 \, \left| \frac{\bar{\lambda}}{2}\, (\bar{\lambda} - \bar{\mu}) + \frac{(\alpha \beta - 1)}{q}\, t^2 - i \, \lambda \, ( \bar{\lambda} - \bar{\mu})\, B  \right|^2 \\
    &= \frac{\alpha^2 \beta^2}{q^2}\, t^4 \, \left( (\bar{\lambda} - \bar{\mu})^2 \, (\frac{\bar{\lambda}^2}{4} + \lambda^2\, B^2) + \frac{\alpha \beta - 1}{q}\, t^2 \, \bar{\lambda}\, (\bar{\lambda} - \bar{\mu}) + \frac{(\alpha \beta - 1)^2}{q^2}\, t^4  \right) \\
    &= \frac{\alpha^2 \beta^2}{q^2}\, t^4 \, \left( (\bar{\lambda} - \bar{\mu})(\frac{t^2}{q}\, (\bar{\lambda} - \bar{\mu}) + \frac{\alpha \beta - 1}{q}\, t^2\, \bar{\lambda}) +  \frac{(\alpha \beta - 1)^2}{q^2}\, t^4 \right) \\
    &= \frac{\alpha^2 \beta^2}{q^3}\, t^6 \, \left((\bar{\lambda} - \bar{\mu})\, (\alpha \beta \, \bar{\lambda} - \bar{\mu}) + \frac{(\alpha \beta - 1)^2}{q}\, t^2 \right) \,.
\end{align*}
Since the denominator is the same for $S_V\,,S_U$ and $S_W\,$, this final form is helpful in all three computations.

We now focus on the numerator, and more precisely on its real part. The previous computation gives us 
$$
D_+^* \, D_-^* = \frac{\alpha \beta}{q}\, t^2 \, \left(\frac{\bar{\lambda}}{2}\, (\bar{\lambda} - \bar{\mu}) + \frac{\alpha \beta - 1}{q}\, t^2 + i\, \lambda \, (\bar{\lambda} - \bar{\mu})\, B\right)\,.
$$
Moreover, we have
\begin{align*}
    N_+\,D_- - N_-\, D_+ &= \alpha\, (A - i\, B) \, (\tilde{A}\, (D_- - D_+) + i \, \tilde{B}\, (D_- + D_+)) \\
    &= - \frac{2 \alpha \beta}{q} \, t \, \tilde{B} \, (A - i \, B) \, \left(\lambda\,B + i \, (\frac{\bar{\lambda}}{2} + \alpha\, t)\right)\,,
\end{align*}
after some simplifications using $\mu\, \tilde{A} - \frac{\bar{\mu}}{2} = \frac{\beta}{q}\, t\,$ . Also,
$$
D_+^*\,D_-^*\,(A - i\,B) = \left(\frac{\lambda - ct}{2q \lambda}\,t\right)\,(\bar{\lambda} - \bar{\mu}) + \frac{\alpha \beta - 1}{q}\, t^2\, A
 + i \, \frac{t}{q}\, (\bar{\lambda} - \bar{\mu})\, B - i\, \frac{\alpha \beta - 1}{q}\, t^2\, B\,,$$
using $\frac{\bar{\lambda}}{2}\, A + \lambda \, B^2 = \frac{\lambda - ct}{2q \lambda}\,t$ and $\lambda\, A - \frac{\bar{\lambda}}{2} = \frac{t}{q}\,$. We can now compute the real part of the entire numerator, which, after some simplifications, is 
$$
\frac{2 \alpha^2 \beta^2}{q^3}\, t^5 \, B \, \tilde{B}\, \left((1 - \alpha)\, \bar{\mu} + \alpha (1 - \beta)\, \bar{\lambda} + (1 - \alpha \beta)\, (\alpha + \frac{1}{q})\, t\right)\,.
$$
We end up with the announced formula,
\begin{align*}
    \bar{V}(\mu, \lambda , t) &= \frac{1}{2\pi^2 \alpha \, \rho(\lambda,t)\, \tilde{\rho}(\mu,t)}\, \Re\left[S_V^+ - S_V^-\right] \\
    &= q\,\frac{(1-\alpha)\, t \,\bar{\mu} + \alpha \, (1 - \beta)\, t \, \bar{\lambda} + (1 - \alpha \beta)\, (\alpha + \frac{1}{q})\, t^2}{(1 - \alpha \beta)^2\, t^2 + q\,(\bar{\lambda} - \bar{\mu})\, (\alpha \beta \, \bar{\lambda} - \bar{\mu})} \,.
\end{align*}

\bibliographystyle{plain}
\bibliography{ref}

\end{document}